
\magnification=1200



\font\bfit=cmbxti10


\font\sc=cmcsc10


\font\sans=cmss10




\font \fivesans               = cmss10 at 5pt

\font \ninei                  = cmmi9

\font \ninesy                 = cmsy9
\font \ninett                 = cmtt9
\font \sevensans              = cmss10 at 7pt

\font \sixi                   = cmmi6
\font \eightrm						= cmr8

\font \sixsy                  = cmsy6

\font \tams                   = cmmib10
\font \tamss                  = cmmib10 scaled 700

\font \tbmss                  = cmmib10 scaled 600

\font \tensans                = cmss10

\font\sixi=cmmi10 at 6pt

\font\sixsy=cmsy10 at 6pt

%

\font \tbmss                  = cmmib10 scaled 600
\newdimen\oldparindent\oldparindent=1.5em
\parindent=1.5em
\skewchar\ninei='177 \skewchar\sixi='177
\skewchar\ninesy='60 \skewchar\sixsy='60
\hyphenchar\ninett=-1
\def\newline{\hfil\break}%
\newfam\sansfam
\textfont\sansfam=\tensans\scriptfont\sansfam=\sevensans
\scriptscriptfont\sansfam=\fivesans
\def\sans{\fam\sansfam\tensans}
\def\bbbA{{\mathchoice {\setbox0=\hbox{$\displaystyle\sans A$}\hbox{\raise
				0.15\ht0\hbox to0pt{\kern0.4\wd0\vrule height0.8\ht0\hss}\box0}}
				{\setbox0=\hbox{$\textstyle\sans A$}\hbox{\raise
				0.15\ht0\hbox to0pt{\kern0.4\wd0\vrule height0.8\ht0\hss}\box0}}
				{\setbox0=\hbox{$\scriptstyle\sans A$}\hbox{\raise
				0.15\ht0\hbox to0pt{\kern0.4\wd0\vrule height0.7\ht0\hss}\box0}}
				{\setbox0=\hbox{$\scriptscriptstyle\sans A$}\hbox{\raise
				0.15\ht0\hbox to0pt{\kern0.4\wd0\vrule height0.7\ht0\hss}\box0}}}}

\def\bbbC{{\mathchoice {\setbox0=\hbox{$\displaystyle\sans C$}\hbox{\raise
				0.15\ht0\hbox to0pt{\kern0.4\wd0\vrule height0.8\ht0\hss}\box0}}
				{\setbox0=\hbox{$\textstyle\sans C$}\hbox{\raise
				0.15\ht0\hbox to0pt{\kern0.4\wd0\vrule height0.8\ht0\hss}\box0}}
				{\setbox0=\hbox{$\scriptstyle\sans C$}\hbox{\raise
				0.15\ht0\hbox to0pt{\kern0.4\wd0\vrule height0.7\ht0\hss}\box0}}
				{\setbox0=\hbox{$\scriptscriptstyle\sans C$}\hbox{\raise
				0.15\ht0\hbox to0pt{\kern0.4\wd0\vrule height0.7\ht0\hss}\box0}}}}

\def\bbbG{{\mathchoice {\setbox0=\hbox{$\displaystyle\sans G$}\hbox{\raise
				0.15\ht0\hbox to0pt{\kern0.4\wd0\vrule height0.8\ht0\hss}\box0}}
				{\setbox0=\hbox{$\textstyle\sans G$}\hbox{\raise
				0.15\ht0\hbox to0pt{\kern0.4\wd0\vrule height0.8\ht0\hss}\box0}}
				{\setbox0=\hbox{$\scriptstyle\sans G$}\hbox{\raise
				0.15\ht0\hbox to0pt{\kern0.4\wd0\vrule height0.7\ht0\hss}\box0}}
				{\setbox0=\hbox{$\scriptscriptstyle\sans G$}\hbox{\raise
				0.15\ht0\hbox to0pt{\kern0.4\wd0\vrule height0.7\ht0\hss}\box0}}}}

\def\bbbJ{{\mathchoice {\setbox0=\hbox{$\displaystyle\sans J$}\hbox{\raise
				0.15\ht0\hbox to0pt{\kern0.4\wd0\vrule height0.8\ht0\hss}\box0}}
				{\setbox0=\hbox{$\textstyle\sans J$}\hbox{\raise
				0.15\ht0\hbox to0pt{\kern0.4\wd0\vrule height0.8\ht0\hss}\box0}}
				{\setbox0=\hbox{$\scriptstyle\sans J$}\hbox{\raise
				0.15\ht0\hbox to0pt{\kern0.4\wd0\vrule height0.7\ht0\hss}\box0}}
				{\setbox0=\hbox{$\scriptscriptstyle\sans J$}\hbox{\raise
				0.15\ht0\hbox to0pt{\kern0.4\wd0\vrule height0.7\ht0\hss}\box0}}}}

\def\bbbO{{\mathchoice {\setbox0=\hbox{$\displaystyle\sans O$}\hbox{\raise
				0.15\ht0\hbox to0pt{\kern0.4\wd0\vrule height0.8\ht0\hss}\box0}}
				{\setbox0=\hbox{$\textstyle\sans O$}\hbox{\raise
				0.15\ht0\hbox to0pt{\kern0.4\wd0\vrule height0.8\ht0\hss}\box0}}
				{\setbox0=\hbox{$\scriptstyle\sans O$}\hbox{\raise
				0.15\ht0\hbox to0pt{\kern0.4\wd0\vrule height0.7\ht0\hss}\box0}}
				{\setbox0=\hbox{$\scriptscriptstyle\sans O$}\hbox{\raise
				0.15\ht0\hbox to0pt{\kern0.4\wd0\vrule height0.7\ht0\hss}\box0}}}}

\def\bbbQ{{\mathchoice {\setbox0=\hbox{$\displaystyle\sans Q$}\hbox{\raise
				0.15\ht0\hbox to0pt{\kern0.4\wd0\vrule height0.8\ht0\hss}\box0}}
				{\setbox0=\hbox{$\textstyle\sans Q$}\hbox{\raise
				0.15\ht0\hbox to0pt{\kern0.4\wd0\vrule height0.8\ht0\hss}\box0}}
				{\setbox0=\hbox{$\scriptstyle\sans Q$}\hbox{\raise
				0.15\ht0\hbox to0pt{\kern0.4\wd0\vrule height0.7\ht0\hss}\box0}}
				{\setbox0=\hbox{$\scriptscriptstyle\sans Q$}\hbox{\raise
				0.15\ht0\hbox to0pt{\kern0.4\wd0\vrule height0.7\ht0\hss}\box0}}}}
\def\bbbR{{\sans I\!R}}
\def\bbbS{{\mathchoice{\setbox0=\hbox{$\displaystyle\sans S$}
				\hbox{\raise0.5\ht0\hbox
				to0pt{\kern0.35\wd0\vrule height0.45\ht0\hss}\hbox
				to0pt{\kern0.55\wd0\vrule height0.5\ht0\hss}\box0}}
				{\setbox0=\hbox{$\textstyle        \sans S$}\hbox{\raise0.5\ht0\hbox
				to0pt{\kern0.35\wd0\vrule height0.45\ht0\hss}\hbox
				to0pt{\kern0.55\wd0\vrule height0.5\ht0\hss}\box0}}
				{\setbox0=\hbox{$\scriptstyle\sans S$}\hbox{\raise0.5\ht0\hbox
				to0pt{\kern0.35\wd0\vrule height0.45\ht0\hss}\raise0.05\ht0\hbox
				to0pt{\kern0.5\wd0\vrule height0.45\ht0\hss}\box0}}
				{\setbox0=\hbox{$\scriptscriptstyle\sans S$}\hbox{\raise0.5\ht0\hbox
				to0pt{\kern0.4\wd0\vrule height0.45\ht0\hss}\raise0.05\ht0\hbox
				to0pt{\kern0.55\wd0\vrule height0.45\ht0\hss}\box0}}}}
\def\bbbT{{\mathchoice {\setbox0=\hbox{$\displaystyle\sans
				T$}\hbox{\hbox to0pt{\kern0.3\wd0\vrule height0.9\ht0\hss}\box0}}
				{\setbox0=\hbox{$\textstyle\sans T$}\hbox{\hbox
				to0pt{\kern0.3\wd0\vrule height0.9\ht0\hss}\box0}}
				{\setbox0=\hbox{$\scriptstyle\sans T$}\hbox{\hbox
				to0pt{\kern0.3\wd0\vrule height0.9\ht0\hss}\box0}}
				{\setbox0=\hbox{$\scriptscriptstyle\sans T$}\hbox{\hbox
				to0pt{\kern0.3\wd0\vrule height0.9\ht0\hss}\box0}}}}
\def\bbbU{{\mathchoice {\setbox0=\hbox{$\displaystyle\sans U$}\hbox{\raise
				0.15\ht0\hbox to0pt{\kern0.4\wd0\vrule height0.8\ht0\hss}\box0}}
				{\setbox0=\hbox{$\textstyle\sans U$}\hbox{\raise
				0.15\ht0\hbox to0pt{\kern0.4\wd0\vrule height0.8\ht0\hss}\box0}}
				{\setbox0=\hbox{$\scriptstyle\sans U$}\hbox{\raise
				0.15\ht0\hbox to0pt{\kern0.4\wd0\vrule height0.7\ht0\hss}\box0}}
				{\setbox0=\hbox{$\scriptscriptstyle\sans U$}\hbox{\raise
				0.15\ht0\hbox to0pt{\kern0.4\wd0\vrule height0.7\ht0\hss}\box0}}}}
\def\bbbV{{\mathchoice {\hbox{$\sans\textstyle V\kern-0.4em V$}}
				{\hbox{$\sans\textstyle V\kern-0.4em V$}}
				{\hbox{$\sans\scriptstyle V\kern-0.3em V$}}
				{\hbox{$\sans\scriptscriptstyle V\kern-0.2em V$}}}}
				\def\qed{\ifmmode\sq\else{\unskip\nobreak\hfil
				\penalty50\hskip1em\null\nobreak\hfil\sq
				\parfillskip=0pt\finalhyphendemerits=0\endgraf}\fi}
\def\bbbW{{\mathchoice {\hbox{$\sans\textstyle W\kern-0.4em W$}}
				{\hbox{$\sans\textstyle W\kern-0.4em W$}}
				{\hbox{$\sans\scriptstyle W\kern-0.3em W$}}
				{\hbox{$\sans\scriptscriptstyle W\kern-0.2em W$}}}}
				\def\qed{\ifmmode\sq\else{\unskip\nobreak\hfil
				\penalty50\hskip1em\null\nobreak\hfil\sq
				\parfillskip=0pt\finalhyphendemerits=0\endgraf}\fi}
\def\bbbX{{\mathchoice {\hbox{$\sans\textstyle X\kern-0.4em X$}}
				{\hbox{$\sans\textstyle X\kern-0.4em X$}}
				{\hbox{$\sans\scriptstyle X\kern-0.3em X$}}
				{\hbox{$\sans\scriptscriptstyle X\kern-0.2em X$}}}}
				\def\qed{\ifmmode\sq\else{\unskip\nobreak\hfil
				\penalty50\hskip1em\null\nobreak\hfil\sq
				\parfillskip=0pt\finalhyphendemerits=0\endgraf}\fi}
\def\bbbY{{\mathchoice {\hbox{$\sans\textstyle Y\kern-0.4em Y$}}
				{\hbox{$\sans\textstyle Y\kern-0.4em Y$}}
				{\hbox{$\sans\scriptstyle Y\kern-0.3em Y$}}
				{\hbox{$\sans\scriptscriptstyle Y\kern-0.2em Y$}}}}
				\def\qed{\ifmmode\sq\else{\unskip\nobreak\hfil
				\penalty50\hskip1em\null\nobreak\hfil\sq
				\parfillskip=0pt\finalhyphendemerits=0\endgraf}\fi}
\def\bbbZ{{\mathchoice {\hbox{$\sans\textstyle Z\kern-0.4em Z$}}
				{\hbox{$\sans\textstyle Z\kern-0.4em Z$}}
				{\hbox{$\sans\scriptstyle Z\kern-0.3em Z$}}
				{\hbox{$\sans\scriptscriptstyle Z\kern-0.2em Z$}}}}
				\def\qed{\ifmmode\sq\else{\unskip\nobreak\hfil
				\penalty50\hskip1em\null\nobreak\hfil\sq
				\parfillskip=0pt\finalhyphendemerits=0\endgraf}\fi}

\def\bbba{{\mathchoice {\setbox0=\hbox{$\displaystyle\sans a$}\hbox{\raise
          0.15\ht0\hbox to0pt{\kern0.4\wd0\vrule height0.8\ht0\hss}\box0}}
          {\setbox0=\hbox{$\textstyle\sans a$}\hbox{\raise
          0.15\ht0\hbox to0pt{\kern0.4\wd0\vrule height0.8\ht0\hss}\box0}}
          {\setbox0=\hbox{$\scriptstyle\sans a$}\hbox{\raise
          0.15\ht0\hbox to0pt{\kern0.4\wd0\vrule height0.7\ht0\hss}\box0}}
          {\setbox0=\hbox{$\scriptscriptstyle\sans a$}\hbox{\raise
          0.15\ht0\hbox to0pt{\kern0.4\wd0\vrule
height0.7\ht0\hss}\box0}}}}
\def\bbbb{{\mathchoice {\setbox0=\hbox{$\displaystyle\sans b$}\hbox{\raise
          0.15\ht0\hbox to0pt{\kern0.4\wd0\vrule height0.8\ht0\hss}\box0}}
          {\setbox0=\hbox{$\textstyle\sans b$}\hbox{\raise
          0.15\ht0\hbox to0pt{\kern0.4\wd0\vrule height0.8\ht0\hss}\box0}}
          {\setbox0=\hbox{$\scriptstyle\sans b$}\hbox{\raise
          0.15\ht0\hbox to0pt{\kern0.4\wd0\vrule height0.7\ht0\hss}\box0}}
          {\setbox0=\hbox{$\scriptscriptstyle\sans b$}\hbox{\raise
          0.15\ht0\hbox to0pt{\kern0.4\wd0\vrule
height0.7\ht0\hss}\box0}}}}
\def\bbbc{{\mathchoice {\setbox0=\hbox{$\displaystyle\sans C$}\hbox{\raise
          0.15\ht0\hbox to0pt{\kern0.4\wd0\vrule height0.8\ht0\hss}\box0}}
          {\setbox0=\hbox{$\textstyle\sans C$}\hbox{\raise
          0.15\ht0\hbox to0pt{\kern0.4\wd0\vrule height0.8\ht0\hss}\box0}}
          {\setbox0=\hbox{$\scriptstyle\sans C$}\hbox{\raise
          0.15\ht0\hbox to0pt{\kern0.4\wd0\vrule height0.7\ht0\hss}\box0}}
          {\setbox0=\hbox{$\scriptscriptstyle\sans C$}\hbox{\raise
          0.15\ht0\hbox to0pt{\kern0.4\wd0\vrule
height0.7\ht0\hss}\box0}}}}
\def\bbbd{{\mathchoice {\setbox0=\hbox{$\displaystyle\sans d$}\hbox{\raise
          0.15\ht0\hbox to0pt{\kern0.4\wd0\vrule height0.8\ht0\hss}\box0}}
          {\setbox0=\hbox{$\textstyle\sans d$}\hbox{\raise
          0.15\ht0\hbox to0pt{\kern0.4\wd0\vrule height0.8\ht0\hss}\box0}}
          {\setbox0=\hbox{$\scriptstyle\sans d$}\hbox{\raise
          0.15\ht0\hbox to0pt{\kern0.4\wd0\vrule height0.7\ht0\hss}\box0}}
          {\setbox0=\hbox{$\scriptscriptstyle\sans d$}\hbox{\raise
          0.15\ht0\hbox to0pt{\kern0.4\wd0\vrule
height0.7\ht0\hss}\box0}}}}
\def\bbbe{{\mathchoice {\setbox0=\hbox{$\displaystyle\sans e$}\hbox{\raise
          0.15\ht0\hbox to0pt{\kern0.4\wd0\vrule height0.8\ht0\hss}\box0}}
          {\setbox0=\hbox{$\textstyle\sans e$}\hbox{\raise
          0.15\ht0\hbox to0pt{\kern0.4\wd0\vrule height0.8\ht0\hss}\box0}}
          {\setbox0=\hbox{$\scriptstyle\sans e$}\hbox{\raise
          0.15\ht0\hbox to0pt{\kern0.4\wd0\vrule height0.7\ht0\hss}\box0}}
          {\setbox0=\hbox{$\scriptscriptstyle\sans e$}\hbox{\raise
          0.15\ht0\hbox to0pt{\kern0.4\wd0\vrule
height0.7\ht0\hss}\box0}}}}
\def\bbbf{{\mathchoice {\setbox0=\hbox{$\displaystyle\sans f$}\hbox{\raise
          0.15\ht0\hbox to0pt{\kern0.4\wd0\vrule height0.8\ht0\hss}\box0}}
          {\setbox0=\hbox{$\textstyle\sans f$}\hbox{\raise
          0.15\ht0\hbox to0pt{\kern0.4\wd0\vrule height0.8\ht0\hss}\box0}}
          {\setbox0=\hbox{$\scriptstyle\sans f$}\hbox{\raise
          0.15\ht0\hbox to0pt{\kern0.4\wd0\vrule height0.7\ht0\hss}\box0}}
          {\setbox0=\hbox{$\scriptscriptstyle\sans f$}\hbox{\raise
          0.15\ht0\hbox to0pt{\kern0.4\wd0\vrule
height0.7\ht0\hss}\box0}}}}
\def\bbbg{{\mathchoice {\setbox0=\hbox{$\displaystyle\sans g$}\hbox{\raise
          0.15\ht0\hbox to0pt{\kern0.4\wd0\vrule height0.8\ht0\hss}\box0}}
          {\setbox0=\hbox{$\textstyle\sans g$}\hbox{\raise
          0.15\ht0\hbox to0pt{\kern0.4\wd0\vrule height0.8\ht0\hss}\box0}}
          {\setbox0=\hbox{$\scriptstyle\sans g$}\hbox{\raise
          0.15\ht0\hbox to0pt{\kern0.4\wd0\vrule height0.7\ht0\hss}\box0}}
          {\setbox0=\hbox{$\scriptscriptstyle\sans g$}\hbox{\raise
          0.15\ht0\hbox to0pt{\kern0.4\wd0\vrule
height0.7\ht0\hss}\box0}}}}
\def\bbbu{{\mathchoice {\setbox0=\hbox{$\displaystyle\sans u$}\hbox{\raise
				0.15\ht0\hbox to0pt{\kern0.4\wd0\vrule height0.8\ht0\hss}\box0}}
				{\setbox0=\hbox{$\textstyle\sans u$}\hbox{\raise
				0.15\ht0\hbox to0pt{\kern0.4\wd0\vrule height0.8\ht0\hss}\box0}}
				{\setbox0=\hbox{$\scriptstyle\sans u$}\hbox{\raise
				0.15\ht0\hbox to0pt{\kern0.4\wd0\vrule height0.7\ht0\hss}\box0}}
				{\setbox0=\hbox{$\scriptscriptstyle\sans u$}\hbox{\raise
				0.15\ht0\hbox to0pt{\kern0.4\wd0\vrule height0.7\ht0\hss}\box0}}}}
\def\bbby{{\mathchoice {\setbox0=\hbox{$\displaystyle\sans y$}\hbox{\raise
          0.15\ht0\hbox to0pt{\kern0.4\wd0\vrule height0.8\ht0\hss}\box0}}
          {\setbox0=\hbox{$\textstyle\sans y$}\hbox{\raise
          0.15\ht0\hbox to0pt{\kern0.4\wd0\vrule height0.8\ht0\hss}\box0}}
          {\setbox0=\hbox{$\scriptstyle\sans y$}\hbox{\raise
          0.15\ht0\hbox to0pt{\kern0.4\wd0\vrule height0.7\ht0\hss}\box0}}
          {\setbox0=\hbox{$\scriptscriptstyle\sans y$}\hbox{\raise
          0.15\ht0\hbox to0pt{\kern0.4\wd0\vrule
height0.7\ht0\hss}\box0}}}}
\def\acknow#1{
 \if N\lasttitle\removelastskip\vskip\baselineskip
     \fi
     \bgroup\petit
     \noindent
     {\it Acknowledgement.\ }%
     \ignorespaces#1\par\egroup
     \ignorespaces}
\def\newcenterenvironment#1#2#3#4{\long\def#1##1##2{\removelastskip
\vskip\baselineskip\centerline{#3#2\if!##1!.\else\ ##1.\fi
\ }{#4\ignorespaces##2}\vskip\baselineskip}}
\newcenterenvironment\centit{\bf}{\bf}{\bf}
\def\newenvironment#1#2#3#4{\long\def#1##1##2{\removelastskip
\vskip\baselineskip\noindent{#3#2\if!##1!.\else\ ##1.\fi
\ }{#4\ignorespaces##2}\vskip\baselineskip}}
\newenvironment\lemma{Lemma}{\bf}{\sl}
\newenvironment\proposition{Proposition}{\bf}{\sl}
\newenvironment\theorem{Theorem}{\bf}{\sl}
\newenvironment\extensionoftheorem{Extension of Theorem}{\bf}{\sl}
\newenvironment\corollary{Corollary}{\bf}{\sl}
\newenvironment\conjecture{Conjecture}{\bf}{\sl}
\newenvironment\claim{Claim}{\bf}{\sl}
\newenvironment\algorithm{Algorithm}{\bf}{\rm}
\newenvironment\example{Example}{\bf}{\rm}
\newenvironment\examples{Examples}{\bf}{\rm}
\newenvironment\exercise{Exercise}{\bf}{\rm}
\newenvironment\problem{Problem}{\bf}{\rm}
\newenvironment\solution{Solution}{\bf}{\rm}
\newenvironment\definition{Definition}{\bf}{\rm}
\newenvironment\remark{Remark}{\bf}{\rm}
\newenvironment\comment{Comment}{\bf}{\rm}
\newenvironment\note{Note}{\bf}{\rm}
\newenvironment\notation{Notation}{\bf}{\rm}
\newenvironment\application{Application}{\bf}{\rm}
\newenvironment\case{Case}{\bf}{\rm}
\newenvironment\question{Question}{\bf}{\rm}
\newenvironment\answer{Answer}{\bf}{\rm}
\newenvironment\algorithm{Algorithm}{\bf}{\rm}
\newenvironment\rule{Rule}{\bf}{\rm}
\newenvironment\result{Result}{\bf}{\rm}
\newenvironment\terminology{Terminology}{\bf}{\rm}
\long\def\proof{\removelastskip\vskip\baselineskip\noindent{\sl
Proof.\quad}\ignorespaces}
\let\lasttitle=N
\def\title{\textfont0=\tenbf \scriptfont0=\sevenbf
\scriptscriptfont0=\fivebf
\textfont1=\tams \scriptfont1=\tamss \scriptscriptfont1=\tbmss}
 \def\centitle#1#2{\if N\lasttitle\else\vskip-24pt
     \fi
     \vskip24pt plus 4pt minus4pt
     \bgroup
\textfont0=\tenbf \scriptfont0=\sevenbf \scriptscriptfont0=\fivebf
\textfont1=\tams \scriptfont1=\tamss \scriptscriptfont1=\tbmss
     \lineskip=0pt
     \pretolerance=10000
     \centerline
     \bf
     \rightskip 0pt plus 6em
     \setbox0=\vbox{\vskip23pt\def\fonote##1{}%
     \if!#1!\ignorespaces#2
     \else\setbox0=\hbox{\ignorespaces#1\unskip.\enspace}%
     \hangindent=\wd0
     \hangafter=1\box0\ignorespaces#2\fi
     \vskip10pt}%
     \dimen0=\pagetotal\advance\dimen0 by-\pageshrink
     \ifdim\dimen0<\pagegoal
     \dimen0=\ht0\advance\dimen0 by\dp0\advance\dimen0 by
     3\normalbaselineskip
     \advance\dimen0 by\pagetotal
     \ifdim\dimen0>\pagegoal\eject\fi\fi
     \
     \if!#1!\ignorespaces#2
     \else\setbox0=\hbox{\ignorespaces#1\unskip.\enspace}%
     \hangindent=\wd0
     \hangafter=1\box0\ignorespaces#2\fi
     \vskip12pt plus4pt minus4pt\egroup
     \nobreak
     \parindent=0pt
     \everypar={\global\parindent=\oldparindent
     \global\let\lasttitle=N\global\everypar={}}%
     \global\let\lasttitle=A%
     \ignorespaces}
 \def\titlea#1#2{\if N\lasttitle\else\vskip-24pt
     \fi
     \vskip24pt plus 4pt minus4pt
     \bgroup
\textfont0=\tenbf \scriptfont0=\sevenbf \scriptscriptfont0=\fivebf
\textfont1=\tams \scriptfont1=\tamss \scriptscriptfont1=\tbmss
     \lineskip=0pt
     \pretolerance=10000
     \noindent
     \bf
     \rightskip 0pt plus 6em
     \setbox0=\vbox{\vskip23pt\def\fonote##1{}%
     \if!#1!\ignorespaces#2
     \else\setbox0=\hbox{\ignorespaces#1\unskip.\enspace}%
     \hangindent=\wd0
     \hangafter=1\box0\ignorespaces#2\fi
     \vskip10pt}%
     \dimen0=\pagetotal\advance\dimen0 by-\pageshrink
     \ifdim\dimen0<\pagegoal
   \dimen0=\ht0\advance\dimen0 by\dp0\advance\dimen0 by
     3\normalbaselineskip
     \advance\dimen0 by\pagetotal
     \ifdim\dimen0>\pagegoal\eject\fi\fi
     \noindent
     \if!#1!\ignorespaces#2
     \else\setbox0=\hbox{\ignorespaces#1\unskip.\enspace}%
     \hangindent=\wd0
     \hangafter=1\box0\ignorespaces#2\fi
     \vskip12pt plus4pt minus4pt\egroup
     \nobreak
     \parindent=0pt
     \everypar={\global\parindent=\oldparindent
     \global\let\lasttitle=N\global\everypar={}}%
     \global\let\lasttitle=A%
     \ignorespaces}
 \def\titleb#1#2{\if N\lasttitle\else\vskip-24pt
     \fi
     \vskip24pt plus 4pt minus4pt
     \bgroup
     \it
     \lineskip=0pt
     \pretolerance=10000
     \noindent
     \rightskip 0pt plus 6em
     \setbox0=\vbox{\vskip23pt\def\fonote##1{}%
     \noindent
     \if!#1!\ignorespaces#2
     \else\setbox0=\hbox{\ignorespaces#1\unskip.\enspace}%
     \box0%
     \ignorespaces#2\fi
     \vskip6pt}%
     \dimen0=\pagetotal\advance\dimen0 by-\pageshrink
     \ifdim\dimen0<\pagegoal
     \dimen0=\ht0\advance\dimen0 by\dp0\advance\dimen0 by
     2\normalbaselineskip
     \advance\dimen0 by\pagetotal
     \ifdim\dimen0>\pagegoal\eject\fi\fi
     \noindent
     \if!#1!\ignorespaces#2
     \else\setbox0=\hbox{\ignorespaces#1\unskip.\enspace}%
     \box0%
     \ignorespaces#2\fi
     \vskip12pt plus4pt minus4pt\egroup
     \nobreak
     \parindent=0pt
     \everypar={\global\parindent=\oldparindent
     \global\let\lasttitle=N\global\everypar={}}%
     \global\let\lasttitle=B%
     \ignorespaces}
 \def\titlec#1{\if N\lasttitle\else\vskip-\baselineskip
     \fi
     \vskip12pt plus 4pt minus 4pt
     \bgroup
     \it
     \noindent
     \ignorespaces#1\unskip.\ \egroup
     \ignorespaces}
 \def\titlecb#1{\if N\lasttitle\else\vskip-\baselineskip
     \fi
     \vskip12pt plus 4pt minus 4pt
     \bgroup
     \bf
     \noindent
     \ignorespaces#1\unskip.\ \egroup
     \ignorespaces}

 \def\titled#1{\if N\lasttitle\removelastskip\vskip\baselineskip
     \fi
     \bgroup
     \noindent
\textfont0=\tenbf \scriptfont0=\sevenbf \scriptscriptfont0=\fivebf
\textfont1=\tams \scriptfont1=\tamss \scriptscriptfont1=\tbmss
     \tenbf
     \ignorespaces#1\unskip.\ \egroup
     \ignorespaces}



\def\ve{\mathop{\varepsilon}\nolimits}

\def\vp{\mathop{\varpi}\nolimits}

\def\vp{\mathop{\varphi}\nolimits}




\def\tildel{\mathop{\widetilde\delta}\nolimits}

\def\tilQ{\mathop{\widetilde Q}\nolimits}


\def\tilv{\mathop{\widetilde v}\nolimits}

\def\tily{\mathop{\widetilde y}\nolimits}
\def\tilz{\mathop{\widetilde z}\nolimits}


\def\baru{\mathop{\overline u}\nolimits}

\def\barw{\mathop{\overline w}\nolimits}
\def\barx{\mathop{\overline x}\nolimits}


\def\calb{\mathop{\cal B}\nolimits}

\def\calg{\mathop{\cal G}\nolimits}
\def\calh{\mathop{\cal H}\nolimits}

\def\calm{\mathop{\cal M}\nolimits}

\def\calq{\mathop{\cal Q}\nolimits}

\def\calz{\mathop{\cal Z}\nolimits}

\def\ra{\mathop{\rightarrow}\nolimits}

\def\linebox{\hbox to .50truein{\hrulefill}}
\def\proofend{\hfill{\vrule height .9ex width .8ex depth -.1ex}
	\vskip\baselineskip}

\def\boxit#1{\vskip\baselineskip{\vbox{\hrule\hbox{\vrule\kern3pt
	\vbox{\kern3pt#1\kern3pt}\kern3pt\vrule}\hrule}}\vskip\baselineskip}

\def\iiitem{\par\indent\indent \hangindent3\parindent\textindent}

\newcount\refnumber
\refnumber=1

\def\refnam#1{\xdef#1{
\item{\rm[\the\refnumber]}
}\global\advance\refnumber by 1}

\newcount\eqnumber
\eqnumber=1

\def\new#1{{\rm \the\eqnumber}
\global\advance\eqnumber by 1}

\tolerance=10000
\overfullrule=0pt
\baselineskip=18pt
\parindent=18pt
\parskip=5pt
\nopagenumbers
\footline={\eightrm 
$\;\;$\hfil\folio\hfil 6/24/94}
\def\span{\mathop{\hbox{\rm{span}}}\limits}
\def\ker{\mathop{\hbox{\rm{ker}}}\nolimits}
\def\Im{\mathop{\hbox{\rm{im}}}\nolimits}
\def\sIm{\mathop{\hbox{\eightrm{im}}}\nolimits}
\def\id{\mathop{\hbox{\rm{id}}}\nolimits}
\def\gla{\hbox{\sans{GLA}}}
\def\dgla{\hbox{\sans{DGLA}}}
\def\opluslim{\mathop{\oplus}\limits}

\centerline{\bf The Spectral Flow of the Odd Signature Operator and Higher
Massey Products}

\vskip\baselineskip

\centerline{\sc Paul A. Kirk and Eric P. Klassen}

\vskip\baselineskip
\vskip\baselineskip

\titlea{0}{Introduction}

Given an analytic path $A_t$ of flat connections on a
principal $U(N)$-bundle $P$ over a closed odd-dimensional
manifold $M^{2\ell-1}$, how can we calculate the spectral
flow of the corresponding path $D_t$ of signature operators on
an associated vector bundle $E\ra M$?  In this paper we give an
algebraic-topological answer to this question in terms of
cohomology, cup products and higher Massey products.  We
summarize here our technique.

 Let $E_{\bbbC}$ denote the complexified adjoint Lie-algebra bundle
associated to  $P$. For the purposes of this summary, we will assume that
$E=E_{\bbbC}$; the more general case will be dealt with in section 7. Let
$d_t:\Omega^p(M;E_{\bbbC})\ra\Omega^{p+1}(M;E_{\bbbC})$  denote the
exterior
derivative corresponding to $A_t$ for  each $t$.  At $t=0$, we wish to
 calculate
the dimension of  $\ker(D_0)$, which gives the number of eigenvalues
$\lambda_{\alpha}(t)$ of $D_t$ passing through $0$ at $t=0$.
Then, for each of these $\lambda_{\alpha}(t)$ which vanish
at $t=0$, we need to calculate the first non-vanishing
derivative of $\lambda_{\alpha}(t)$ at $t=0$.  Because the
analyticity of $A_t$ implies that each $\lambda_{\alpha}(t)$
is analytic, this information will give a complete
description of the spectral flow of $D_t$ near $t=0$.

To obtain this information, we set up a sequence of cochain
complexes $\{\calg^*_n,\delta_n\}$, for $n=0,1,2,\ldots$.
Note that $\calg^i_n$ denotes the $i^{\hbox{\eightrm{th}}}$
cochain group of the $n^{\hbox{\eightrm{th}}}$ cochain complex,
and the coboundary operator of this complex is $\delta_n$.  $\calg_0^*$ is
just
the deRham complex of $M$ with respect to $d_0$, i.e.,
$\calg^i_0=\Omega^i(M;E_{\bbbC})$ and $\delta_0=d_0$.  For
each $n$, the chain groups of $\calg^*_{n+1}$ are the
cohomology groups of $\calg^*_n$, i.e.,
$\calg^i_{n+1}=H^i(\calg^*_n,\delta_n)$.  For example,
$\calg^i_1=H^i(M;E_{\bbbC})_{A(0)}$.

The coboundary operator $\delta_n$ is defined by
$\delta_n[u]=(-
1)^{n+1}[Q_{\calm}(u,{\underbrace{a_1,\ldots,a_1}_n})]$.

In this formula, $Q_{\calm}$ denotes an $(n+1)$-fold higher
Massey product of the cocycle $u$ with $n$ copies of the
$1$-cocycle $a_1$, where $a_1$ is the derivative ${d\over
dt}\big|_{t=0}A(t)$ of the path of flat connections.  The
Massey system $\calm$ which we use to compute this product
is required to be $A(t)$-compatible (see
section 6 for the definition).

As a special case, $\delta_1$ is just the map
$H^i(M;E_{\bbbC})_{A(0)}\ra H^{i+1}(M;E_{\bbbC})_{A(0)}$
given by the ``cup product'' $u\mapsto[u,a_1]$.  Thus at this
first level, no Massey system is needed to compute the boundary operator.

The significance of these cochain complexes
$(\calg^*_n,\delta_n)$ is as follows.  The number of
eigenvalues $\lambda_{\alpha}(t)$ of $D_t$ which vanish to
order $n-1$ at $t=0$ is given by
$\sum\limits_{k{\hbox{\eightrm{
even}}}}\dim(\calg^k_{n})$.  Of those eigenvalues which
vanish to order $n-1$, the difference between the number whose
$n^{\hbox{\eightrm{th}}}$ derivative is positive and the
number whose $n^{\hbox{\eightrm{th}}}$ derivative is
negative is the signature of the Hermitian form
$$Q_n:\calg^{\ell-1}_n\times\calg^{\ell-1}_n\ra H^{2\ell-
1}(M;\bbbC)=\bbbC$$
given by
$$Q_n(v,w)=\cases{
i(-1)^{\ell+1\over 2}\delta_nv\cdot w&if $\ell$
is odd\cr
(-1)^{\ell\over 2}\delta_nv\cdot w&if $\ell$ is
even\cr}$$
The ``$\cdot$'' in this expression denotes the cup product
constructed using the natural Hermitian structure on
$u(N)\otimes\bbbC$.

In essence, what this says is that if you understand the cup
product and higher Massey product structures on
$H^*(M;E_{\bbbC})_{A(0)}$, then you understand the local
behavior of those eigenvalues of the signature operator
$D_t$ which pass through $0$ at $t=0$.

As one might expect, if you know the path $A(t)$ to
$n^{\hbox{\eightrm{th}}}$ order, then you know the behavior of
the eigenvalues to $n^{\hbox{\eightrm{th}}}$ order.  For
example, to understand the eigenvalues' behavior to first
order, you need only know the derivative $a_1$ of the path
$A(t)$.  To construct the Massey system needed to compute
$\delta_n$, you need to know $A(t)$ to order $n$ (again, see
section 6 for details).

The first order form $Q_1$ was first defined in [KK1] for 3-manifolds and
in
[KK2] for odd-dimensional manifolds.  In [FL], Farber and Levine gave a
definition for the higher-order forms $\{Q_n\}$ involving a linking
form on the
cohomology of $M$ with coefficients in a power series ring. Our present
definition of the $Q_n$ has the advantage that it describes these forms in
terms of cup products and higher Massey products without referring directly
 to
power series (though power series are used in the proofs).

Our sequence $\{\calg^*_n\}$ of chain complexes corresponds to the spectral
sequence which Farber obtains in [Fa] using an analytic deformation of an
elliptic complex.  In our exposition, instead of using the technology of
spectral sequences, we define this sequence of chain complexes ``by hand''
in
order to show clearly that the differentials $\delta_n$ can be described as
higher Massey products.

The paper is organized as follows:

In section 1, we outline our strategy for calculating
derivatives of eigenvalues of an analytic path $D_t$ of self-adjoint
operators
on a Hilbert space $V$, which involves  defining a nested sequence of
subspaces
$V=V_0\supseteq  V_1\supseteq V_2\supseteq\ldots$ and Hermitian forms $B_n$
on $V_n$.

In section 2, given an analytic path of operators
$d_t:V\rightarrow V$ such that $d_t\circ d_t=0$ for all $t$,
we construct a sequence of subquotients $\calg_n$ of $V$ and
maps $\delta_n:\calg_n\rightarrow\calg_n$ such that
$\delta_n\circ\delta_n=0$,
$\calg_{n+1}={\ker\delta_n\over\sIm\delta_n}$, and
$\delta_n$ depends on the $n^{\hbox{\eightrm{th}}}$ order
expansion of $d_t$.

In section 3, using an inner product on $V$ in addition to
the path $d_t:V\ra V$ with $d_t\circ d_t=0$, we construct a
nested sequence of subspaces $\calh_n$ of $V$ such that
$\calh_n\cong\calg_n$, and maps
$\tildel_n:\calh_n\ra\calh_n$ corresponding to the maps
$\delta_n:\calg_n\ra\calg_n$, such that
$\calh_n=\ker\tildel_n\cap\ker\tildel^*_n$.  This gives
``harmonic representatives'' of the elements of $\calg_n$.

In section 4, we show that for the operator $T_t=d_t+d^*_t$,
the spaces $V_n$ constructed in section 1 are exactly the
same as the spaces $\calh_n$ constructed in section 3.

In section 5, for a principal $U(N)$-bundle $P\ra M^{2\ell-
1}$, a connection $A$ on $P$, and the adjoint
$u(N)\otimes\bbbC$-bundle $E_{\bbbC}\ra M$, we define the
signature operator
$D:\Omega^{\hbox{\eightrm{even}}}(M;E_{\bbbC})\ra
\Omega^{\hbox{\eightrm{even}}}(M;E_{\bbbC})$.  Under the
assumption that $A_t$ is a path of flat connections and
$D_t$ is the corresponding path of operators, we show that
$V_n(D)=\opluslim\limits_{k{\hbox{\eightrm{
even}}}}\calh^k_n$, where $V_n(D)$ are the subspaces
corresponding to $D_t$ as defined in section 1 and
$\calh_n=\oplus\calh^i_n$ are the ``harmonic'' subspaces
defined in section 3.  We also write down a formula for the
forms $B_n(D)$ on $V_n(D)$ (as defined in section 1) in
terms of the operator $\tildel_n$ and $*$.  (The dimension
of $V_n(D)$ is the number of eigenvalues of $D_t$ which
vanish to order $n-1$; the $n^{\hbox{\eightrm{th}}}$
derivatives of these eigenvalues are given by the
eigenvalues of $B_n(D)$.)  Finally we define a ``reduced
version'' $Q_n$ of $B_n(D)$.  $Q_n$ has domain $\calh^{\ell-
1}_n=\calg^{\ell-1}_n$, is defined only in terms of $\delta_n$ and cup
products (hence is independent of the metric on $M$), and has the same
signature
as $B_n(D)$.

In section 6, we show that the operators $\delta_n$ are
actually just higher Massey products (in the sense of
Retakh) for the differential graded Lie algebra
$\Omega^*(M;E_{\bbbC})$.

Finally, in section 7, we indicate  how to handle the case of the odd
signature
operator on a general  bundle $E\ra M$.  This involves ``hybrid''  Massey
products taking into account the action of the Lie algebra bundle
$E_{\bbbC}$ on
$E$.

\titlea{1}{Derivatives of vanishing eigenvalues of a path of self-adjoint
operators}

Let $D_t:V\ra V$ be an analytic path of closed, self-adjoint operators on
the complex Hilbert   space $V$. We will assume that each $D_t$ has compact
resolvent, so that the spectrum of $D_t$ is a discrete subset of $\bbbR$
and
each eigenvalue has finite multiplicity.

We recall from [Ka] that the path $D_t$ is called {\it analytic} if there
exists a Hilbert space $W$ and a path $Q_t$ in $Bd(W,V)$ so that $Q_t$ is
analytic (with respect to the operator norm topology), the domain of $D_t$
is
the range of $Q_t$, and the composite $D_t\circ Q_t$ is an analytic path in
$Bd(W,V)$.

\note{}{The way this set-up will occur later in the paper is that $V$
 will be a
space of sections of a Hermitian vector bundle over a closed
Riemannian manifold
and  $D_t$ will be an analytic path of first order elliptic operators which
are self-adjoint with respect to the $L^2$ inner product. For details, see
section 5. Of course, for the simplest example one may think of $V$ as
finite
dimensional and $D_t$ as an analytic path of Hermitian matrices.}

  By [Ka],
there is an analytic family $\{\vp_{\alpha}(t)\}_{\alpha\in A}$ of
orthonormal Hilbert bases for $V$ and a collection of analytic
real valued functions $\{\lambda_{\alpha}(t)\}_{\alpha\in A}$ such
that for all $t\in(-\ve,\ve)$,
$$D_t\vp_{\alpha}(t)=\lambda_{\alpha}(t)\vp_{\alpha}(t).$$
We wish to calculate the first non-vanishing derivative of
$\lambda_{\alpha}(t)$ at $t=0$ for those $\lambda_{\alpha}$'s
satisfying $\lambda_{\alpha}(0)=0$.  We accomplish this by
defining and diagonalizing a sequence of bilinear forms on a
nested sequence of subspaces of $V$.

First, define the subsets $A_0,A_1,A_2,\ldots$ of $A$ by
$$\openup3\jot\eqalign{
A_0&=A\cr
A_n&=\left\{\alpha\in A:{d^i\over
dt^i}\Big\vert_{t=0}\;\lambda_{\alpha}(t)=0{\hbox{\rm{ for all
}}}i<n\right\}\cr}$$
Define subspaces $V_n\subseteq V$ by
$$V_n=\span\limits_{\alpha\in A_n}\{\vp_{\alpha}(0)\}$$
For each $n$, define a Hermitian bilinear form $B_n:V_n\times
V_n\ra\bbbC$ as follows.

Given $v,w\in V_n$, write $v=\sum\limits_{\alpha\in
A_n}c_{\alpha}\vp_{\alpha}(0)$.  Define
$$B_n(v,w)=\left\langle{d^n\over
dt^n}\bigg\vert_{t=0}\;D_t\;\sum\limits_{\alpha\in
A_n}\;c_{\alpha}\vp_{\alpha}(t),w\right\rangle.$$
{}From this definition it follows immediately that
$$B_n(\vp_{\alpha}(0),\vp_{\beta}(0))=\lambda^{(n)}_{\alpha}(0)
\delta_{\alpha\beta}.$$
Hence $B_n$ is Hermitian and its eigenvalues are the
$n^{\hbox{\eightrm{th}}}$ derivatives of those eigenvalues of $D_t$
which vanish to order $n-1$.  It also follows that $V_{n+1}$ is
the subspace on which $B_n$ is degenerate.

\titlea{2}{A sequence of chain complexes}

Let $V$ be a vector space and let $d_t:V\ra V$ be an analytic path of
linear transformations such that $d_t\circ d_t=0$ for all $t$ near
$0$. To avoid any technicalities, we will assume that
$d_t=d_0+\sum_{i=1}^\infty A_it^i$ where $\sum_{i=1}^\infty A_it^i$ is a
convergent power series in $Bd(V,V)$. and $d_0$ is some closed linear
operator.  The domain of $d_t$ is taken to be the same for all $t$.  (In
the main example to keep in mind,
$V$ will be the $L^2$ closure of $\oplus_nV_n$  where
$(V_n,d_t)$ is an elliptic complex over a closed smooth manifold.  In
particular, $\ker d_0/\Im d_0$ will be finite dimensional.)

 We will define inductively a sequence
$\calg_0,\calg_1,\calg_2,\ldots$ of vector spaces and
$\delta_n:\calg_n\ra\calg_n$ of linear transformations as follows.

First, let $\calg_0=V$ and $\delta_0=d_0$.

Inductively, let $\calg_n=\ker(\delta_{n-1})/\Im(\delta_{n-1})$.

To define $\delta_n$, note that $\calg_n$ is a subquotient of $V$.
Hence an element of $\calg_n$ may be written as $[u]$ where $u\in
V$.  Define
$$\delta_n([u])=\left[{d^n\over
dt^n}\bigg\vert_{t=0}\;d_t(u+u_1t+\ldots+u_{n-1}t^{n-1})\right]$$
where $u_1,\ldots,u_{n-1}\in V$ are chosen so that for all $i<n$,
$${d^i\over dt^i}\bigg\vert_{t=0}\;d_t(u+u_1t+\ldots+u_{n-1}t^{n-
1})=0.$$

\theorem{2.1}{The spaces $\calg_n$ and operators
$\delta_n:\calg_n\ra\calg_n$ are well-defined and
$\delta_n\circ\delta_n=0$ for all $n$.}

The rest of the section is devoted to proving this theorem.

\proof Since $d_t$ is analytic we may write
$d_t=d+A_1t+A_2t^2+\ldots$ where each $A_i$ is linear and $d=d_0$.  In this
notation, if we apply $d_t$ to a formal power series in $t$ with
coefficients in $V$, we obtain another such formal power series.
With this in mind, we define a sequence $\{T_n\}^{\infty}_{n=0}$
of linear relations on $V$ as follows:  $T_0$ is the function
$d_0$.  Define $uT_nv$ if and only if there exist
$u_1,\ldots,u_{n-1}\in V$ such that $d_t(u+u_1t+\ldots+u_{n-
1}t^{n-1})=vt^n+w(t)t^{n+1}$ for some $w(t)\in V[[t]]$.  As usual,
we define
$$\openup3\jot\eqalign{
{\hbox{\rm{Domain}}}(T_n)&=\{u\in V:uT_nv{\hbox{\rm{ for some
}}}v\in V\}\cr
{\hbox{\rm{Range}}}(T_n)&=\{v\in V:uT_nv{\hbox{\rm{ for some
}}}u\in V\}\cr}$$
We now state four lemmas concerning these relations.

\lemma{2.2}{Suppose $n\geq 1$, $k\geq 0$ and $u,v\in V$.  Then
$uT_nv\Rightarrow vT_k0$.}

\lemma{2.3}{Suppose $n\geq 1$.  Then $v\in{\hbox{\rm{
Range}}}(T_n)\Leftrightarrow 0T_{n+1}v$.}

\lemma{2.4}{Suppose $v=du$.  Then for all $n\geq 0$, there exists
$w_n\in V$ such that $vT_n\,dw_n$.}

\lemma{2.5}{

\item{(a)}$u\in{\hbox{\rm{ Domain}}}(T_1)\Leftrightarrow du=0$.

\item{(b)}If $n\geq 2$, then $u\in{\hbox{\rm{
Domain}}}(T_n)\Leftrightarrow uT_{n-1}\,dw$ for some $w\in V$.

}

\noindent{\sl Proofs.}

\item{(2.2)}$uT_nv$ implies we may choose $u_1,\ldots,u_{n-1}\in V$
such that $d_t(u+u_1t+\ldots+u_{n-1}t^{n-
1})=vt^n+w(t)t^{n+1}=(v+w(t)t)t^n$.  Since $d^2_t=0$, it follows
that $d_t(v+w(t)t)=0$ which implies that $vT_k0$ for all $k\geq
0$.

\itemitem{(2.3)$(\Rightarrow)$}Since $uT_nv$, there exist
$u_1,\ldots,u_{n-1}$ such that $d_t(u+u_1t+\ldots+u_{n-1}t^{n-
1})=vt^n+w(t)t^{n+1}$.  Multiplying both sides by $t$ immediately
implies $0T_{n+1}v$.

\itemitem{$(\Leftarrow)$}Just reverse the above argument.

\item{(2.4)}Define $v_i=A_iu$ for $i\geq 1$, so
$d_tu=(v+v_1t+v_2t^2+\ldots)$.  Since $d^2_t=0$,
$d_t(v+v_1t+v_2t^2+\ldots)=0$ which implies
$d_t(v+v_1t+\ldots+v_{n-1}t^{n-1})=-d_t(v_nt^n+\ldots)$ which we
can write
$$d_t(v+v_1t+\ldots+v_{n-1}t^{n-1})=d(-v_n)t^n+w(t)t^{n+1}.$$
Hence we may let $w_n=-v_n$ and the lemma follows.

\itemitem{(2.5)  (a)}There exists $v\in V$ such that $uT_1v$.

\iiitem{$\Leftrightarrow$}There exists $v\in V$ such that
$d_tu=vt+w(t)t^2$.

\iiitem{$\Leftrightarrow$}$du=0$.

\iiitem{(b)$(\Rightarrow)$}$u\in{\hbox{\rm{ Domain}}}(T_n)$
implies there exist $u_1,\ldots,u_{n-1}$ such that
$$d_t(u+u_1t+\ldots+u_{n-1}t^{n-1})=w(t)t^n$$
for some $w(t)\in
V[[t]].$ It follows that
$$\openup3\jot\eqalign{
d_t(u+u_1t+\ldots+u_{n-2}t^{n-2})
&=-d_tu_{n-1}t^{n-1}+y(t)t^n\cr
&=d(-u_{n-1})t^{n-1}+\tily(t)t^n\cr}$$
where $\tily(t)\in V[[t]]$.  So $uT_{n-1}\,d(-u_{n-1})$.

\iiitem{$(\Leftarrow)$}Suppose $uT_{n-1}\,dw$.  Then
$$d_t(u+u_1t+\ldots+u_{n-2}t^{n-2})=dwt^{n-1}+y(t)t^n.$$
Let $u_{n-1}=-w$.  Then
$$\openup3\jot\eqalign{
d_t(u+u_1t+\ldots+u_{n-1}t^{n-1})
&=-d_twt^{n-1}+dwt^{n-1}+y(t)t^n\cr
&=\tily(t)t^n.\cr}$$
So $u\in{\hbox{\rm{ Domain}}}(T_n)$.

Having proven these lemmas, we make the following definitions:

\definition{}{Let $\calz_n$ and $\calb_n$ denote the following
subspaces of $V$ for $n=1,2,3,\ldots$.
$$\openup3\jot\displaylines{
\calz_1=\,\ker(d)\cr
\calb_1=\,\Im(d).\cr}$$
For $n\geq 2$,
$$\openup3\jot\eqalign{
\calz_n&=\{u\in V:uT_{n-1}\,dw{\hbox{\rm{ for some }}}w\}\cr
\calb_n&={\hbox{\rm{ Range}}}(T_{n-1})+dV.\cr}$$

}

\proposition{2.6}{For all $n\geq 1$,

\item{(a)}$\calb_n\subseteq\calz_n$.

\item{(b)}$\calz_{n+1}\subseteq\calz_n$.

\item{(c)}$\calb_n\subseteq\calb_{n+1}$.

}

This proposition follows easily from the above lemmas and we omit
the proof.

Define $G_n=\calz_n/\calb_n$.

\proposition{2.7}{$T_n$ induces a well-defined linear operator
$\tau_n:G_n\ra
G_n$ and $\tau_n\circ \tau_n=0$.}

\proof $(n=1)$.  $u\in\calz_1\Rightarrow du=0\Rightarrow
u\in{\hbox{\rm{ Domain}}}(T_1)$ by Lemma 2.5.  If $u\in\calz_1$ and
$uT_1v$ then, by Lemma 2.2, $v\in\calz_1$.

To see that $T_1$ induces a well-defined $\tau_1:G_1\ra G_1$, we note
first that there are no choices involved in defining $T_1$, and
then note that by Lemma 2.4, $T_1(\calb_1)\subseteq\calb_1$.

$(n\geq 2)$.  If $u\in\calz_n$ then, by Lemma 2.5, $u\in{\hbox{\rm{
Domain}}}(T_n)$.  By Lemma 2.2, if $uT_nv$ then $v\in\calz_n$.

If $u\in\calb_n$, then we can write $u=u_1+u_2$ where there
exists an $x$ such that $xT_{n-1}u_1$ and there exists a $y$ such
that $dy=u_2$.  It follows by Lemmas 2.2 and 2.4 that $uT_nv$ where
$v\in\calb_n$.  If it is also true that $uT_nv\rq$, then $0T_n(v-
v\rq)$.  Hence $v-v\rq\in{\hbox{\rm{ Range}}}(T_{n-1})$ (by Lemma 2.3)
so $v\rq$ is also in $\calb_n$.  Hence $\tau_n:G_n\ra G_n$ is
well-defined.

Lemma 2.2 implies that $\tau_n\circ \tau_n=0$ for all $n$.

\proposition{2.8}{$G_{n+1}=\,\ker\,\tau_n/\Im\,\tau_n$.}

\proof $\tau_n:G_n\ra G_n$ where $G_n=\calz_n/\calb_n$.  It is true
that $\ker\,\tau_n=\calz_{n+1}/\calb_n$ and
$\Im(\tau_n)=\calb_{n+1}/\calb_n$.  We will show that:
$\ker\,\tau_n\subseteq \calz_{n+1}/\calb_n$, which is the only
non-trivial part.

Suppose $u\in\calz_n$ and $\tau_n[u]=0$.  Then $uT_nv$ where
$v\in\calb_n$.  Write $v=v_1+v_2$ where
$v_1\in{\hbox{\rm{ Range }}}T_{n-1}$ and $v_2\in dV$.  Since
$v_1\in{\hbox{\rm{ Range }}}T_{n-1}$, Lemma 2.3 implies that
$0T_nv_1$, so $uT_nv_2$.  Hence $u\in\calz_{n+1}$, and
$\ker\,\tau_n\subseteq\calz_{n+1}/\calb_n$.\proofend

Letting $\calg_n=G_n$ and $\delta_n=n!\tau_n$, we have now proven
theorem 2.1, since $\delta_n$ and $\calg_n$ satisfy the formula
and inductive definition given at the beginning of this section.

\vskip4ex

Some of the results of this section can also be described in terms
of a spectral sequence; this is essentially the content of Theorem
6.1 of [Fa].  We give a brief alternative description here for
completeness.  (Note that for much of what follows in our paper,
 it will still
be important to have the explicit description of $\tau_n$ which we
 have given
above.)

Let $M=V[[t]]$ denote the formal power series with coefficients in
$V$, and filter $M$ by $F^p=t^pM$.  The differential $d_t:M\rightarrow M$
preserves the filtration  and thus one obtains a spectral
sequence.  In terms of an exact couple one takes $D_1=\oplus
D_1^p$ where $D^p_1= H^*(F^p;d_t)$ and $E_1=\oplus E_1^p$ where
$E^p_1=H^*(F^{p-1}/F^{p}; d_t)$.  Then
$$\matrix{
D_1 &          & \rightarrow &         & D_1\cr
    &          &             &         &    \cr
    & \nwarrow &             &\swarrow &    \cr
    &          &             &         &    \cr
    &          &     E_1     &         &    \cr}$$
is an exact couple.
   One identifies the higher
$E_r$ by taking $E_r^p=Z_r^p/B_r^p$ where
$$Z^p_r= k^{-1}( \hbox{Im} i^{r-1}:D_1^{p+r}\to D^{p+1}_1)  $$
and
$$B^p_r= j( \ker i^{r-1}:D_1^{p}\to D^{p-r+1}_1).$$
(See [Mc].) It is then a routine exercise to show inductively that
$G_n\cong E^p_n$ for any $p\ge n-1$, and that $\tau_n:G_n\to G_n$
coincides with the induced differential $d_n:E^p_n\to E^{p-n}_n$ for $p\ge
2n-1$.

\titlea{3}{Hodge Theory}

Assume we are in the situation of the last chapter; i.e.,
$d_t:V\ra V$ is an analytic path of linear transformation such
that $d^2_t=0$.  Assume in addition that $V$ has a Hermitian inner product
$\langle\quad,\quad\rangle$.  For each $t$, let $d^*_t:V\ra V$
denote the adjoint of $d_t$.  Clearly $d^*_t\circ d^*_t=0$.

We assume furthermore that $d_t^*$ is a closed operator with the same
domain as
$d_t$.   Recall that $d_t=d_0 + A(t)$ where $A(t)$ is an analytic path in
$Bd(V,V)$.  Thus $d_t^*=d_0^*+A(t)^*$ is an analytic path.

Recall from section 2 that
$\calz_1=\,\ker\,d$ and $\calz_n=\{u\in V:uT_{n-1}\,dw{\hbox{\rm{
for some }}}w\in V\}$ for $n\geq 2$.  By Lemma 2.5 in section 2,
$\calz_n={\hbox{\rm{ Domain }}}T_n$.  Hence an alternative
description of $\calz_n$ (for $n\geq 1$) is
$$\openup3\jot\eqalign{
&\{u\in V:\;\exists\;u_1,\ldots,u_{n-1}\in V{\hbox{\rm{ such that }}}\cr
&\qquad{d^k\over dt^k}\bigg\vert_{t=0}d_t(u+tu_1+\ldots+t^{n-
1}u_{n-1})=0{\hbox{\rm{ for all }}}k<n\}\cr
&=\{u\in V:\;\exists\;u_1,\ldots,u_{n-1}\in V{\hbox{\rm{ such that }}}\cr
&\qquad d_t(u+tu_1+\ldots+t^{n-1}u_{n-1})=v(t)t^n{\hbox{\rm{ for some
}}}v(t)\in V[[t]]\}.\cr}$$
Define $\calz^*_n\subseteq V$ analogously using $d^*_t$:
$$\openup3\jot\eqalign{
\calz^*_1&=\,\ker\,d^*\cr
\calz^*_n&=\{u\in V:\;\exists\;u_1,\ldots,u_{n-1}\in V{\hbox{\rm{ such that
}}}\cr
&\quad d^*_t(u+tu_1+\ldots+t^{n-1}u_{n-
1})=v(t)t^n{\hbox{\rm{ for some }}}v(t)\in V[[t]]\}.\cr}$$
Define $\calh_n=\calz_n\cap\calz^*_n$.  Note that
$\calz_1\supseteq\calz_2\supseteq\calz_3\supseteq\ldots$ and
$\calz^*_1\supseteq\calz^*_2\supseteq\calz^*_3\supseteq\ldots$ are
nested sequences, hence so is
$\calh_1\supseteq\calh_2\supseteq\calh_3\supseteq\ldots$.

\lemma{3.1}{$\calz^*_n\subseteq\calb^{\perp}_n$.}

\proof Let $u\in\calz^*_n$; i.e., assume that there exist
$u_1,\ldots,u_{n-1}\in V$ such that $d^*_t(u+tu_1+\ldots+t^{n-
1}u_{n-1})=q(t)t^n$ for some $q(t)\in V[[t]]$.  We need to show
$u\perp dV$ and $u\perp{\hbox{\rm{ Range}}}(T_{n-1})$.  Since
$u\in\calz^*_n\subseteq\calz^*_1=\,\ker\,d^*$ it follows that
$u\perp dV$.  Now let $v\in{\hbox{\rm{ Range}}}(T_{n-1})$; i.e.,
assume there exist $w,w_1,\ldots,w_{n-2}$ such that
$$d_t(w+tw_1+\ldots+t^{n-2}w_{n-2})=vt^{n-1}+z(t)t^n.$$
Clearly
$$\openup3\jot\eqalign{
\langle u,v\rangle t^{n-1}+(\quad)t^n
&=\langle u+tu_1+\ldots+t^{n-1}u_{n-1},vt^{n-1}+z(t)t^n\rangle\cr
&=\langle u+tu_1+\ldots+t^{n-1}u_{n-1},d_t(w+tw_1+\ldots+t^{n-
2}w_{n-2})\rangle\cr
&=\langle q(t)t^n,w+tw_1+\ldots+t^{n-2}w_{n-2}\rangle\cr
&=\langle q(t),w+tw_1+\ldots,t^{n-2}w_{n-2}\rangle
t^n.\cr}$$
\proofend

So $\langle u,v\rangle =0$.  Since
$\calz^*_n\subseteq\calb^{\perp}_n$, it follows that the map
$$\Phi_n:\calz_n\cap\calz^*_n\ra{\calz_n\over\calb_n}$$
is injective.  Recall from section 2 that $\calz_n={\hbox{\rm{
Domain}}}(T_n)$ and ${\hbox{\rm{Range}}}(T_n)\subseteq\calz_n$.
Also the indeterminancy of $T_n$ is contained in $\calb_n$.  It
follows that we have a well-defined homomorphism
$\tildel_n:\calz_n\cap\calz^*_n\ra\calz_n\cap\calz^*_n$
defined by $\tildel_n=pr_{\calh_n}\circ T_n$.  Explicitly,
$$\tildel_n(u)=pr_{\calh_n}\left({d^n\over dt^n}\bigg\vert_{t=0}d_t
(u+tu_1+\ldots,t^{n-1}u_{n-1})\right)$$
where $u_1,\ldots,u_{n-1}$ are chosen so that
$${d^k\over dt^k}\bigg\vert_{t=0}d_t
(u+tu_1+\ldots,t^{n-1}u_{n-1})=0$$
for all $k<n$.  By symmetry, the homomorphism
$\tildel^*_n:\calz_n\cap\calz^*_n\ra\calz_n\cap\calz^*_n$ defined
by
$$\tildel^*_nu=pr_{\calh_n}{d^n\over
dt^n}\bigg\vert_{t=0}d^*_t(u+tu_1+\ldots+t^{n-1}u_{n-1})=0,$$
where $u_1,\ldots,u_{n-1}$ are chosen so that
$${d^k\over dt^k}\bigg\vert_{t=0}d^*_t(u+tu_1+\ldots+t^{n-1}u_{n-1})=0$$
for all $k<n$, is also well-defined.

We now justify the notation ``$\tildel^*_n$.''

\proposition{3.2}{$\tildel_n$ and $\tildel^*_n$ are adjoint on
$\calh_n$.}

\proof Let $u,v\in\calh_n$.  Choose $u_1,\ldots,u_{n-1}$ so that
$${d^k\over dt^k}\bigg\vert_{t=0}d_t(u+tu_1+\ldots+t^{n-1}u_{n-
1})=0$$
for all $k<n$ and choose $v_1,\ldots,v_{n-1}$ so that
$${d^k\over dt^k}\bigg\vert_{t=0}d^*_t(v+tu_1+\ldots+t^{n-1}v_{n-
1})=0$$
for all $k<n$.   Note that
$$\openup3\jot\displaylines{
\langle
d_t(u+tu_1+\ldots+t^{n-1}u_{n-1}),(v+tu_1+\ldots+t^{n-1}v_{n-1})\rangle\cr
-\langle (u+tu_1+\ldots+t^{n-1}u_{n-1}),d^*_t(v+tu_1+\ldots+t^{n-
1}v_{n-1})\rangle=0.\cr}$$
Differentiating this expression $n$ times at $t=0$ yields
$$
\left\langle {d^n\over
dt^n}\bigg\vert_{t=0}d_t(u+tu_1+\ldots+t^{n-1}u_{n-1}),v\right\rangle
=\left\langle u,{d^n\over
dt^n}\bigg\vert_{t=0}d^*_t(v+tv_1+\ldots+t^{n-1}v_{n-
1})\right\rangle.$$
Since $u$, $v\in\calh_n$, inserting $pr_{\calh_n}$ in the
appropriate places doesn't alter the values of these inner
products.  Hence
$$\langle\tildel_nu,v\rangle=\langle
u,\tildel^*_nv\rangle.$$
\proofend

\theorem{3.3}{For all $n\geq 1$,
$$\Phi_n:\calz_n\cap\calz^*_n\ra{\calz_n\over\calb_n}=\calg_n$$
is an isomorphism and $\Phi_n\delta_n\Phi^{-1}_n=\tildel_n$.}

\proof We prove the theorem by induction on $n$.

First consider the case $n=1$,  In this case
$\calz^*_1=\,\ker\,d^*=\,\Im(d)^{\perp}=\calb^{\perp}_1$, so
$$\Phi_1:\calz_1\cap\calz^*_1\ra{\calz_1\over\calb_1}$$
is an isomorphism.  The equation $\Phi_1\delta_1\Phi^{-
1}_1=\tildel_1$ follows immediately from the definitions.

For the inductive step, assume we've proven that
$\Phi_n:\calz_n\cap\calz^*_n\ra\calz_n/\calb_n$ is an isomorphism
and $\Phi_n\delta_n\Phi^{-1}_n=\tildel_n$.  Note that (by section
2),
$$\openup3\jot\eqalign{
{\calz_{n+1}\over\calb_{n+1}}={\ker\,\delta_n\over\Im\,\delta_n}
&\cong{\ker\,\tildel_n\over\Im\,\tildel_n}\cr
&\cong\,\ker\,\tildel_n\cap\,\ker\,\tildel^*_n.\cr}$$
(The last isomorphism follows because $\tildel_n\circ\tildel_n=0$
since $\tildel_n=\Phi_n\circ\delta_n\circ\Phi^{-1}_n$.)

All of these isomorphisms are induced by inclusion.

\lemma{}{$\ker(\tildel_n)\cap\,\ker(\tildel^*_n)=\calz_{n+1}\cap
\calz^*_{n+1}$.}

\proof We begin by showing that
$\calz_{n+1}\cap\calh_n\subseteq\,\ker(\tildel_n)$.  If
$u\in\calz_{n+1}\cap\calh_n$, we can choose $u_1,\ldots,u_n$ so that
$d_t(u+tu_1+\ldots+t^nu_n)=z(t)t^{n+1}$.  Hence
$$\openup3\jot\eqalign{
d_t(u+tu_1+\ldots+t^{n-1}u_{n-1})
&=-d_tu_nt^n+z(t)t^{n+1}\cr
&=-du_nt^n+\tilz(t)t^{n+1}.\cr}$$
So $\tildel_n(u)=pr_{\calh_n}(-n!\,du_n)=0$ since
$du_n\in\calb_n\subseteq(\calz^*_n)^{\perp}$.

By symmetry we also conclude that
$\calz^*_{n+1}\cap\calh_n\subseteq\,\ker(\tildel^*_n)$; hence
$$\calz_{n+1}\cap\calz^*_{n+1}=(\calz_{n+1}\cap\calh_n)\cap(\calz^
*_{n+1}\cap\calh_n)\subseteq\,\ker(\tildel_n)\cap\,\ker(\tildel^*_
n).$$
We now show that $\ker(\tildel_n)\subseteq\calz_{n+1}$.

If $u\in\,\ker(\tildel_n)$, we can choose $u_1,\ldots,u_{n-1}$ so
that
$$d_t(u+tu_1+\ldots+t^{n-1}u_{n-1})=vt^n+(\quad)t^{n+1}
\eqno{(*)}$$
where $pr_{\calh_n}(v)=0$.  Since
$\calz^*_n\subseteq\calb^{\perp}_n$ and by the inductive
assumption that
$$\Phi_n:\calz_n\cap\calz^*_n\ra{\calz_n\over\calb_n}$$
is an isomorphism, it follows that $v\in\calb_n$.  (The inductive
assumption implies that $\calz_n$ is
the orthogonal direct sum of $\calh_n$ and $\calb_n$.)  Hence
$v=\tilv+dw$ where $\tilv\in{\hbox{\rm{ Range }}}T_{n-1}$.  By
Lemma 2.3 in section 2, $\tilv\in {\hbox{\rm{ Range }}}T_{n-
1}\Rightarrow 0T_n\tilv$.  Since (*) above can be written $uT_nv$,
it follows by the linearity of the relation that $uT_n\,dw$ which
means $u\in\calz_{n+1}$, proving that
$\ker(\tildel_n)\subseteq\calz_{n+1}$.  By symmetry we also
conclude that $\ker(\tildel^*_n)\subset\calz^*_{n+1}$, hence
$\ker(\tildel_n)\cap\,\ker(\tildel^*_n)\subseteq\calz_{n+1}\cap
\calz^*_{n+1}$, completing the proof that
$\ker(\tildel_n)\cap\,\ker(\tildel^*_n)=\calz_{n+1}\cap\calz^*_{n+
1}$.\proofend

Since we've already shown that
$${\calz_{n+1}\over\calb_{n+1}}\cong\,\ker\,\tildel_n\cap\,\ker\,
\tildel^*_n,$$
the proof that
$$\calh_{n+1}=\calz_{n+1}\cap\calz^*_{n+1}\cong{\calz_{n+1}\over
\calb_{n+1}}$$
is now complete.  It immediately follows that there is an
orthogonal decomposition
$\calz_{n+1}=\calh_{n+1}\oplus\calb_{n+1}$, and from the
definitions of $\delta_n,\tildel_n$ we obtain
$$\Phi_{n+1}\delta_{n+1}\Phi^{-1}_{n+1}=\tildel_{n+1},$$
completing the induction step and the theorem.\proofend

\titlea{4}{The operator $d+d^*$}

We continue with the assumptions of the previous section. Thus $V$ is a
Hilbert space,
$d_t:V\ra V$ is an analytic path of linear operators, and
$d_t\circ d_t=0$ for all $t$ near $0$.   The adjoint $d_t^*$ has
the same domain as $d_t$ (which has the same domain as $d_0$).

 Define
$T_t=d_t+d^*_t$ with domain equal to the domain of $d_t$.  This is
clearly an analytic path of self-adjoint operators. For the remainder of
this section, we also assume that for each $t$ the operator $T_t$ has
compact  resolvent. (This is true, for example, if $d_t$ is a path of
exterior derivatives on the de Rham complex of a compact manifold since in
this case $d_t+d_t^*$ is an elliptic operator.) It follows that we may
apply the apparatus of section 1 to the operator $T_t$ to obtain subspaces
$V_n$. The goal of  this section is to show that the subspaces $V_n$ (as
defined in section 1) coincide with the  subspaces
$\calh_n$ (as defined in section 3).

Recall from section 1 that we have an analytic path
$\{\vp_{\alpha}(t)\}_{\alpha\in A}$ of orthonormal bases of $V$
and analytic paths $\{\lambda_{\alpha}(t)\}$ of real eigenvalues
satisfying
$$T_t\vp_{\alpha}(t)=\lambda_{\alpha}(t)\vp_{\alpha}(t).$$
Recall that for $n=0,1,2,\ldots$
$$A_n=\left\{\alpha\in A:{d^k\over
dt^k}\bigg\vert_{t=0}\lambda_{\alpha}(t)=0{\hbox{\rm{ for all
}}}k<n\right\}$$
and
$$V_n=\,\span\limits_{\alpha\in A_n}\{\vp_{\alpha}(0)\}.$$

\theorem{4.1}{$V_n=\calh_n$.}

\lemma{4.2}{Let $u(t)$ and $v(t)$ be two smooth paths in $V$ such
that for all $t$, $\langle u(t),v(t)\rangle=0$.  Let $n\geq 0$.
Then
$${d^k\over dt^k}\bigg\vert_{t=0}(u(t)+v(t))=0$$
for all $k\leq n$ if and only if
$${d^k\over dt^k}\bigg\vert_{t=0}u(t)=0\quad{\hbox{\rm{ and }}}\quad
{d^k\over dt^k}\bigg\vert_{t=0}v(t)=0$$
for all $k\leq n$.}

\noindent{\sl Proof of Lemma.}  The direction $(\Leftarrow)$ is
obvious.  We prove $(\Rightarrow)$ by induction.

For $n=0$, $u(0)+v(0)=0$ implies $\langle u(0),u(0)\rangle=\langle
u(0),u(0)+v(0)\rangle=0$, so $u(0)=0$ and similarly $v(0)=0$.

Now assume the lemma is true for $n$, and suppose
$${d^k\over dt^k}\bigg\vert_{t=0}(u(t)+v(t))=0$$
for all $k\leq n+1$.  By the induction assumption,
$${d^k\over dt^k}\bigg\vert_{t=0}u(t)=0\quad{\hbox{\rm{ and }}}\quad
{d^k\over dt^k}\bigg\vert_{t=0}v(t)=0$$
for all $k\leq n$.

Differentiate the equation
$$\langle u(t),v(t)\rangle=0$$
$2n+2$ times at $t=0$.  After
cancelling all terms which vanish by the induction assumption, we
are left with
$$\pmatrix{2n+2\cr  \noalign{\vskip5pt}  n+1\cr}
\left\langle {d^{n+1}\over
dt^{n+1}}\bigg\vert_{t=0}u(t),
{d^{n+1}\over dt^{n+1}}\bigg\vert_{t=0}v(t)\right\rangle=0.$$
Since
$${d^{n+1}\over dt^{n+1}}\bigg\vert_{t=0}u(t)+
{d^{n+1}\over dt^{n+1}}\bigg\vert_{t=0}v(t)=0,$$
it follows that each of these terms are $0$, finishing the proof
of the lemma.\proofend

\vskip\baselineskip

\noindent{\sl Proof of Theorem 4.1.}  First we show
$V_n\subseteq\calh_n$.  Let $v\in V_n$ and write
$v=\sum\limits_{\alpha\in A_n} c_{\alpha}\vp_{\alpha}(0)$.  Let
$v(t)=\sum\limits_{\alpha\in A_n} c_{\alpha}\vp_{\alpha}(t)$.  By
definition of $A_n$,
$${d^k\over dt^k}\bigg\vert_{t=0}T_tv(t)=0$$
for all $k<n$.  Since $T_tv(t)=d_tv(t)+d^*_tv(t)$ and
$\Im(d_t)\perp\,\Im(d^*_t)$ for all $t$, the previous lemma
implies that
$${d^k\over dt^k}\bigg\vert_{t=0}d_tv(t)=0\quad{\hbox{\rm{ and }}}\quad
{d^k\over dt^k}\bigg\vert_{t=0}d^*_tv(t)=0$$
for all $k< n$.  Hence $v\in\calz_n\cap\calz^*_n=\calh_n$.

We will now show by induction that $\calh_n\subseteq V_n$ for all
$n\geq 1$.  For $n=1$, clearly
$\calh_1=\,\ker\,d\cap\,\ker\,d^*=\,\ker(d+d^*)=V_1$.  Assume
that $\calh_n=V_n$; i.e., $\calz_n\cap\calz^*_n=V_n$.  Let
$v\in\calh_{n+1}$.  Since $\calh_{n+1}\subseteq\calh_n=V_n$, we can
write $v=\sum\limits_{\alpha\in A_n}c_{\alpha}\vp_{\alpha}(0)$.
Define $v(t)=\sum\limits_{\alpha\in
A_n}c_{\alpha}\vp_{\alpha}(t)$.  As in the first part of this
proof, we may conclude by Lemma 4.2 that
$${d^k\over dt^k}\bigg\vert_{t=0}d_tv(t)=0\quad{\hbox{\rm{ and }}}\quad
{d^k\over dt^k}\bigg\vert_{t=0}d^*_tv(t)=0$$
for all $k< n$.  Hence we may use (the degree $n-1$ truncation of)
$v(t)$ to compute $\tildel_n(v)$ and $\tildel^*_n(v)$.  Keep in
mind that since $v\in\calh_{n+1}$,
$\tildel_n(v)=\tildel^*_n(v)=0$.  Let $\beta\in A_n-A_{n+1}$;
i.e., $\lambda^{(k)}_{\beta}(0)=0$ for $k<n$ but
$\lambda^{(n)}_{\beta}(0)\not=0$.  Compute
$$\openup3\jot\eqalign{
\langle v,\lambda^{(n)}_{\beta}(0)\vp_{\beta}(0)\rangle
&=\left\langle v,{d^n\over
dt^n}\bigg\vert_{t=0}T_t\vp_{\beta}(t)\right\rangle\cr
&={d^n\over dt^n}\bigg\vert_{t=0}\langle
v(t),T_t\vp_{\beta}(t)\rangle\cr
&={d^n\over dt^n}\bigg\vert_{t=0}\langle
T_tv(t),\vp_{\beta}(t)\rangle\cr
&={d^n\over dt^n}\bigg\vert_{t=0}\langle
d_tv(t)+d^*_tv(t),\vp_{\beta}(t)\rangle\cr
&=\langle\tildel_n(v)+\tildel^*_n(v),\vp_{\beta}(0)\rangle\cr
&=0\cr}$$
Hence for all $\beta\in A_n-A_{n+1}$, $\langle
v,\vp_{\beta}(0)\rangle=0$ and we conclude that $v\in V_{n+1}$,
completing the proof of the theorem.\proofend

\note{}{(1)  A more direct way to prove this theorem
would be to prove simply that $\calz^*_n=\calb^{\perp}_n$.  We were able
to do this for $n=1$ and $2$, but we are not sure if it is true
for higher $n$.

\noindent(2)  Since for each $n$, $V_n=\calh_n\cong\calg_n$, we have shown
that $\dim(V_n)$ does not depend on
the choice of inner product, but only on the path $d_t$.}

\titlea{5}{The odd signature operator}

In this section we recall the definition of the signature operator
on an odd-dimensional manifold and apply the material of the
previous section to obtain information about the derivatives of
its eigenvalues.

Let $M$ be a closed oriented Riemannian manifold of dimension
$2\ell-1$, $P\ra M$ a principal $U(N)$-bundle and $E\ra M$ the
$u(N)$-bundle associated to $P$ by the adjoint representation
$U(N)\ra O(u(N))$.  Denote by $E_{\bbbC}$ the complexification
$E\otimes\bbbC$.  If $A$ is a connection on $P$, let
$$d_A:\Omega^p(M;E_{\bbbC})\ra\Omega^{p+1}(M;E_{\bbbC})$$
denote the corresponding exterior derivative.  The standard inner
product on $u(N)$ extends to a Hermitian inner product on
$u(N)\otimes\bbbC$ given by $\langle a\otimes z,b\otimes
w\rangle=-tr(ab)\otimes z\barw$.  We use this inner product to
define a ``dot product'' of forms
$$\Omega^p(M;E_{\bbbC})\otimes\Omega^q(M;E_{\bbbC})\ra\Omega^{p+q}
(M;\bbbC)$$
which we denote by $\alpha\otimes\beta\mapsto\alpha\cdot\beta$,
and which is obtained by wedging the ``form part'' and applying
the inner product $\langle\;\;,\;\;\rangle$ to the coefficients.
Analogously, the Lie bracket operation on $u(N)\otimes\bbbC$
(defined by $[a\otimes z,b\otimes w]=[a,b]\otimes zw$) gives rise
to a Lie bracket of forms
$$\Omega^p(M;E_{\bbbC})\otimes\Omega^q(M;E_{\bbbC})\ra\Omega^{p+q}
(M;E_{\bbbC}),$$
which we denote by $\alpha\otimes\beta\mapsto[\alpha,\beta]$.

In order to make use of the Riemannian structure (which we haven't
used so far), introduce the Hodge star operator
$$*:\Omega^p(M;E_{\bbbC})\ra\Omega^{2\ell-1-p}(M;E_{\bbbC})$$
and define a Hermitian inner product on $\Omega^p(M;E_{\bbbC})$ by
$$\langle\alpha,\beta\rangle=\int\limits_M\alpha\cdot*\beta\in
\bbbC.$$
The operator $*$ has the following three properties:

\item{(1)}$*$ is an isometry

\item{(2)}$*$ is self-adjoint

\item{(3)}$(*)^2=\id_{\Omega^{*}(M;E_{\bbbC})}$

\noindent Properties (2) and (3) above use the fact that $M$ is
odd-dimensional.  Define the operator
$d^{*}_A:\Omega^p(M;E_{\bbbC})\ra\Omega^{p-1}(M;E_{\bbbC})$ by
$d^{*}_Aw=(-1)^p* d_A* w$.  It is easy to verify that
$\langle
d_A\alpha,\beta\rangle=\langle\alpha,d^{*}_A\beta\rangle$.
Denote by $T_A$ the operator $d_A+d^{*}_A$ on $\Omega^{*}(M;E_{\bbbC})$.
It is well-known that
$T_A:\Omega^{*}(M;E_{\bbbC})\ra\Omega^{*}(M;E_{\bbbC})$ is
a self-adjoint elliptic operator whose domain is the image in
$L^2(\Omega^{*}(M;E_{\bbbC}))$ of the Sobolev space
$L^2_1(\Omega^{*}(M;E_{\bbbC}))$.

To obtain an operator which shares the three properties of $*$
mentioned above, but that also commutes with $T_A$, we define
$\mu:\Omega^p(M;E_{\bbbC})\ra\Omega^{2\ell-1-p}(M;E_{\bbbC})$ by
$\mu(w)=i^{\ell}(-1)^{p(p+1)\over 2}* w$.  Then $\mu$ is a
self-adjoint isometry, $\mu^2=\id_{\Omega^{*}}$ and
$T_A\circ\mu=\mu\circ T_A$ on $\Omega^{*}(M;E_{\bbbC})$.  It
follows that we have an orthogonal direct sum
$\Omega^{*}(M;E_{\bbbC})=\Omega^+(M;E_{\bbbC})\oplus\Omega^-
(M;E_{\bbbC})$ where $\Omega^{\pm}(M;E_{\bbbC})$ is the $\pm1$-
eigenspace of $\mu$.  If we take $L^2$-completions, this
decomposition is an orthogonal direct sum of Hilbert spaces.
Since $T_A$ commutes with $\mu$, $T_A$ preserves this
decomposition.  Define functions
$$\Psi_{\pm}:\Omega^{\hbox{\eightrm{even}}}(M;E_{\bbbC})\ra\Omega^
{\pm}(M;E_{\bbbC})$$
by $\Psi_{\pm}(w)={1\over\sqrt{2}}(w\pm\mu(w))$.  Clearly
$\Psi_{\pm}$ is an isometry of
$\Omega^{\hbox{\eightrm{even}}}(M;E_{\bbbC})$ onto
$\Omega^{\pm}(M;E_{\bbbC})$.  The operator
$D_A:\Omega^{\hbox{\eightrm{even}}}(M;E_{\bbbC})\ra
\Omega^{\hbox{\eightrm{even}}}(M;E_{\bbbC})$ given by
$D_A=\Psi^{-1}_+\circ
T_A\circ\Psi_+$ is called the {\bfit signature operator}
 for the manifold
$M$.  If $w\in\Omega^{2k}(M;E_{\bbbC})$, the explicit formula for
$D_A$ is
$$D_Aw=i^{\ell}(-1)^{k-1}(* d_A-d_A*)w.$$

We now consider the path of operators $D_{A_t}$, where $A_t$ is an
analytic path of flat connections.  To be more precise, we define
a path $A_t$, for $-\ve<t<\ve$, of connections to be analytic if
for all $-\ve<t<\ve$, and for all
$w\in\Omega^{*}(M;E_{\bbbC})$,
$$d_{A_t}w=d_{A_0}w+[a(t),w]$$
where $a(t)=a_0+ta_1+t^2a_2+\ldots$ is a power series in
$\Omega^1(M;E)$ which converges in the $C^r$-norm for all $r$.  Of
course the requirement that $A_t$ be flat can be phrased
$d_{A_t}\circ d_{A_t}=0$ for all $t\in(-\ve,\ve)$.  To simplify
notation, we will write $d_t$, $T_t$ and $D_t$ instead of
$d_{A_t}$, $T_{A_t}$ and $D_{A_t}$.

The material of sections 1--4
applies to the deRham complex
$$d_t:\Omega^{*}(M;E_{\bbbC})\ra\Omega^{*}(M;E_{\bbbC})$$
and the corresponding operators $T_t=d_t+d^{*}_t$ to give the subspaces
$V_n=\calh_n$, the forms $B_n$, the operators $\tildel_n$ and
$\tildel^{*}_n:\calh_n\ra\calh_n$, etc.  We first observe that
these spaces and maps respect the grading of
$\Omega^{*}(M;E_{\bbbC})$.

\lemma{5.1}{For each $n$, the space $\calh_n$ breaks up as an
orthogonal direct sum
$$\calh_n=\opluslim\limits^{2\ell-1}_{i=0}\calh^i_n$$
where $\calh^i_n=\calh_n\cap\Omega^i(M;E_{\bbbC})$.  Also,
$\tildel_n(\calh^i_n)\subseteq\calh^{i+1}_n$ and
$\tildel^{*}_n(\calh^i_n)\subseteq\calh^{i-1}_n$.}

\proof We use induction on $n$.  The lemma is clearly true for
$\calh_0=\Omega^{*}(M;E_{\bbbC})$, $\tildel_0=d_0$ and
$\tildel^{*}_0=d^{*}_0$.

Now assume the lemma is true for $n$.  Because $\tildel_n$ has
degree $1$ and $\tildel^{*}_n$ has degree $-1$, we
know $\ker(\tildel_n)=\opluslim\limits^{2\ell-
1}_{i=1}\,\ker\tildel_n\cap\Omega^i(M;E_{\bbbC})$ and
$\ker(\tildel^{*}_n)=\opluslim\limits^{2\ell-
1}_{i=1}\ker\tildel^{*}_n\cap\Omega^i(M;E_{\bbbC})$.  Since
$\calh_{n+1}=\ker(\tildel_n)\cap\ker(\tildel^{*}_n)$, it
follows that $\calh_{n+1}=\opluslim\limits^{2\ell-1}_{i=0}\calh_{n+1}$.
To see that
$\tildel_{n+1}(\calh^i_{n+1})\subseteq\calh^{i+1}_{n+1}$, recall
that to define $\tildel_{n+1}(u)$, where $u\in\calh^i_{n+1}$, we
choose $u_1,\ldots,u_n$ such that ${d^k\over
dt^k}|_{t=0}d_t(u+u_1t+\ldots+u_nt^n)=0$ for all $k<n$.  This
equation will still hold if we replace each $u_j$ by its
projection in $\Omega^i(M;E_{\bbbC})$ since $d_t$ raises dimension
by $1$.  It follows that $\tildel_{n+1}(u)\in\calh^{i+1}_{n+1}$,
and similarly that $\tildel^{*}_{n+1}(u)\in\calh^{i-
1}_{n+1}$\proofend

In section 4, we showed that $V_n=\calh_n=\opluslim\limits^{2\ell-
1}_{i=0}\calh^i_n$.

We now turn our attention to the path of operators
$$D_t:\Omega^{\hbox{\eightrm{even}}}(M;E_{\bbbC})\ra
\Omega^{\hbox{\eightrm{even}}}(M;E_{\bbbC}).$$
We denote the corresponding
spaces and forms constructed in section 1 by $V_n(D)$ and $B_n(D)$
(to distinguish them from the ones just considered which
correspond to the path of operators $T_t$).

Using the fact that $d^{*}_t=(-1)^p* d_t*$ it is easy
to see that $*$ and $\mu$ both restrict to isometries between
$\calh^i_n$ and $\calh^{2\ell-1-i}_n$ for all $i$, $n$ and that for
$w\in\calh^p_n$, $\tildel^{*}_nw=(-1)^p*\tildel_n* w$.
Using the formula $D_t=\Psi^{-1}\circ T_t\circ\Psi$, we deduce that
$$V_n(D)=\Psi^{-1}(V_n)=\bigoplus^{\ell-1}_{k=0}\calh^{2k}_n.$$

\note{}{The spaces $\calg_n$ and the functions $\delta_n$ defined
in section 2 also respect the grading; i.e.,
$\calg_n=\opluslim\limits^{2\ell-1}_{i=0}\calg^i_n$ and
$\delta_n(\calg^i_n)\subseteq\calg^{i+1}_n$.  Since these spaces
and functions were defined without using the metric on $M$, the
spaces $\calh^i_n$ (hence $V_n$ and $V_n(D)$) are also independent
of the metric (up to isomorphism) since they are canonically
isomorphic to $\calg^i_n$.  Since $V_{n+1}(D)$ is the subspace of
$V_n(D)$ on which $B_n(D)$ is degenerate, it follows that the rank
of $B_n(D)$ is independent of the metric.  In what follows, we
will show that the signature of $B_n(D)$ is also independent of
the metric.}

The following proposition gives a formula for the form $B_n(D)$ on
$V_n(D)$ in terms of $\tildel_n$ and the Hodge star.

\proposition{5.2}{If $v\in\calh^{2k}_n$ and
$w\in\opluslim\limits^{\ell-1}_{i=0}\calh^{2i}_n$ then
$$B_n(D)(v,w)=i^{\ell}(-1)^{k-1}\langle(*\tildel_n-
\tildel_n*)v,w\rangle.$$

}

\proof In the notation of section 1, write
$$v=\sum\limits_{\alpha\in A_n}c_{\alpha}\vp_{\alpha}(0)$$
and let
$v(t)=\sum c_{\alpha}\vp_{\alpha}(t)$.  (Recall that
$V_n(D)=\span\limits_{\alpha\in A_n}\{\vp_{\alpha}(0)\}$.)  From
section 1, $B_n(D)(v,w)={d^n\over dt^n}\big\vert_{t=0}\langle
D_tv(t),w\rangle$.  By definition of $A_n$, ${d^k\over
dt^k}\big\vert_{t=0}D_tv(t)=0$ for all $k<n$.  Since $\Im(*
d_t)\perp\Im(d_t)$ for all $t$, it follows by Lemma 4.2 that
for all $k<n$, ${d^k\over dt^k}\big\vert_{t=0}d_tv(t)=0={d^k\over
dt^k}\big\vert_{t=0}d_t* v(t)$.  Therefore
$\tildel_n(v)=pr_{\calh_n}{d^n\over
dt^n}\big\vert_{t=0}d_tv(t)$ and $\tildel_n(* v)=pr_{\calh_n}{d^n\over
dt^n}\big\vert_{t=0}d_t* v(t)$.  These equations and the definition of
$B_n(D)$ imply that $B_n(D)(v,w)=i^{\ell}(-1)^{k-
1}\langle(*\tildel_n-\tildel_n*)v,w\rangle$.\proofend

The essence of this proposition is:  To
calculate the $n^{\hbox{\eightrm{th}}}$ derivatives of
those eigenvalues of
$D_t$
which vanish (at $t=0$) to order $n-1$, we can simply diagonalize
the signature operator on the even part of the chain complex
$(\calh_n,\tildel_n)$.  The kernel of this new signature operator
is isomorphic to the even part of the cohomology of
$(\calh_n,\tildel_n)$, which is just $\opluslim\limits^{\ell-
1}_{i=0}\calh^{2i}_{n+1}$.

Before showing that the signature of the forms $B_n(D)$ are
invariant of the metric on $M$, we make some further observations
on the path of operators
$d_t:\Omega^{*}(M;E_{\bbbC})\ra\Omega^{*}(M;E_{\bbbC})$.
In section 2, we defined the subspaces $\calz_n$ and $\calb_n$
arising from this type of setup.  It is easy to see that these
subspaces respect the grading; hence we write
$$\calz_n=\opluslim\limits^{2\ell-
1}_{i=1}\calz^i_n\quad{\hbox{\rm{and}}}\quad\calb_n
=\opluslim\limits^{2\ell-1}_{i=0}\calb^i_n$$
where $\calz^i_n=\calz_n\cap\Omega^i(M;E_{\bbbC})$ and
$\calb^i_n=(\calb_n\cap\Omega^i(M;E_{\bbbC})$.  Also
$\calg_n=\opluslim\limits^{2\ell-1}_{i=0}\calg^i_n$ where
$\calg^i_n=\calz^i_n/\calb^i_n$.

\proposition{5.3}{Consider the pairing
$\Omega^k(M;E_{\bbbC})\times\Omega^{2\ell-1-
k}(M;E_{\bbbC})\ra\bbbC$ defined by
$(\alpha,\beta)\mapsto\int\limits_M\alpha\cdot\beta$, which is
$\bbbC$-linear in the first factor and conjugate-linear in the
second.

This pairing induces a well-defined pairing
$$\calg^k_n\times\calg^{2\ell-1-k}_n\ra\bbbC\quad{\hbox{\rm{for
all }}}n,k.$$

}

\proof Since
$\int\limits_M\alpha\cdot\beta={\overline{\int\limits_M\beta
\cdot\alpha}}$, it will suffice to show that if $\alpha\in\calb^k_n$ and
$\beta\in\calz^{2\ell-1-k}$, then
$\int\limits_M\alpha\cdot\beta=0$.  This reduces to showing that
if $\alpha\in{\hbox{\rm{Range}}}(d_0)$ or
$\alpha\in{\hbox{\rm{Range}}}(T_{n-1})$ then
$\int\alpha\cdot\beta=0$.  If $\alpha\in{\hbox{\rm{Range}}}(d_0)$,
i.e., if $\alpha=d_0\gamma$, then $\int_M\alpha\cdot\beta=\int
d_0\gamma\cdot\beta=\pm\int\gamma\cdot d_0\beta=0$ since
$\calz_n\subseteq\ker(d_0)$.  If $\alpha\in{\hbox{\rm{Range}}}(T_{n-
1})$ then we may find a series $\gamma(t)$ in
$\Omega^{k-1}(M;E_{\bbbC})$ such that ${d^i\over
dt^i}\big\vert_{t=0}d_t\gamma(t)=0$ for all $i<n-1$, and ${d^{n-1}\over
dt^{n-1}}\big\vert_{t=0}d_t\gamma(t)=\alpha$.  Since
 $\beta\in\calz^{2\ell-
1-k}_n$, we may find a series $\beta(t)$ in $\Omega^{2\ell-1-
k}(M;E_{\bbbC})$ such that $\beta(0)=\beta$ and ${d^i\over
dt^i}\big\vert_{t=0}d_t\beta(t)=0$ for all $i<n$.  We then compute
$$\openup3\jot\eqalign{\int\limits_M\alpha\cdot\beta
&={d^{n-1}\over dt^{n-
1}}\big\vert_{t=0}\int\limits_Md_t\gamma(t)\cdot\beta(t)\cr
&=\pm{d^{n-1}\over dt^{n-1}}\big\vert_{t=0}\int\limits_M\gamma(t)\cdot
d_t\beta(t)\cr
&=0,\cr}$$
completing the proof.\proofend

We now define a Hermitian form $Q_n:\calg^{\ell-
1}_n\times\calg^{\ell-1}_n\ra\bbbC$ by
$$Q_n(v,w)=\cases{
i(-1)^{\ell+1\over 2}\int\limits_M\delta_nv\cdot w&if $\ell$ is
odd\cr
(-1)^{\ell\over 2}\int\limits_M\delta_nv\cdot w&if $\ell$ is
even\cr}$$
(Recall $\delta_n:\calg^p_n\ra\calg^{p+1}_n$ was defined in
section 2.)  Note that $Q_n$ is defined without reference to the
metric on $M$.

\theorem{5.4}{The forms $B_n(D)$ and $Q_n$ have the same signatures.}

\proof First, consider the case in which $\ell$ is odd.  Recall
$V_n(D)=\opluslim\limits^{\ell-1}_{k=0}\calh^{2k}_n$,  and that
$\{\calh^i_n\}_{0\leq i\leq 2\ell-1}$ is a cochain complex with
boundary $\tildel_n$.  Define the operator $D_{(n)}:V_n(D)\ra
V_n(D)$ by $D_{(n)}v=i^{\ell}(-1)^{k-1}(*\tildel_n-
\tildel_n*)v$, so $B_n(D)(v,w)=\langle D_{(n)}v,w\rangle$.
Write $V_n(D)=X\oplus Y\oplus W$ where
$$\openup3\jot\eqalign{
X&=\left(\opluslim\limits_{0\leq k\leq{\ell-3\over
2}}\calh^{2k}_n\right)\oplus\Im(\tildel_n:\calh^{\ell-
2}_n\ra\calh^{\ell-1}_n)\cr
Y&=\opluslim\limits_{0\leq k\leq{\ell-3\over
2}}\calh^{\ell+1+2k}\cr
W&=\ker(\tildel_n*:\calh^{\ell-1}_n\ra\calh^{\ell+1}_n)\cr}$$
Conceptually, $X=$ all even forms in $\calh_n$ of dimension less
than $\ell-1$ plus the ``exact'' forms of dimension $\ell-1$, $Y=$
all even forms of dimension greater than $\ell-1$ and $W=$ the
``coclosed'' $\ell-1$ forms.  (The words ``exact'' and
``coclosed'' refer to the operator $\tildel_n$.)

We omit the verifications of the following facts, which are
routine:  $D_{(n)}(X)\subseteq Y$, $D_{(n)}(Y)\subseteq X$,
$(D_{(n)}|X)^{*}=(D_{(n)}|Y)$, and $D_{(n)}(W)\subseteq W$.
It follows immediately that $B_n(D)$ has the same signature as
$B_n(D)|W$.  For $v,w\in W$, note that
$$\openup3\jot\displaylines{
B_n(D)(v,w)=\langle D_{(n)}v,w\rangle=i^{\ell}(-
1)^{\ell}\langle*\tildel_nv,w\rangle\cr
-i^{\ell}\langle\tildel_nv,* w\rangle=i(-1)^{\ell+1\over
2}\int\tildel_nv\cdot w=Q_n(v,w),\cr}$$
where in the last step we are thinking of $v$ and $w$ as
representing classes in $\calg_n$.

Since $\Im(\tildel_n:\calh^{\ell-2}_n\ra\calh^{\ell-1}_n)$ is
contained in the degeneracy of $Q_n$, it follows that $Q_n$,
$B_n(D)|W$, and $B_n(D)$ all have the same signature, completing
the proof for $\ell$ odd.

If $\ell$ is even the method of proof is the same, only the space
$W$ is replaced by the space
$$U=\ker(\tildel_n:\calh^{\ell}_n\ra\calh^{\ell+1}_n).$$

It follows that $B_n(D)$ has the same signature as the form
$(v,w)\mapsto(-1)^{\ell\over 2}\int\tildel_n* v\cdot * w$ or $U$
which, since $*$ is an isometry, has the same signature as
$Q_n$ on $\calg^{\ell-1}_u$.  This completes the proof of
Theorem 5.4.\proofend

\titlea{6}{Massey products}

Thus far we have demonstrated the existence of a sequence of
chain complexes $\{(\calg^i_n,\delta_n)\}$ such that
$(\calg^i_0,\delta_0)$ is just the deRham complex
$(\Omega^i(M;E_{\bbbC}),d_0)$ and such that the chain groups
of the complex $(\calg^i_n,\delta_n)$ are the homology
groups of the previous complex $(\calg^i_{n-1},\delta_{n-
1})$.  Furthermore, the rank of
$\calg^{\hbox{\eightrm{even}}}_n=\opluslim\limits^{\ell-
1}_{k=0}\calg^{2k}_n$ is the number of eigenvalues of the
signature operator $D_t$ which vanish to order $n-1$ at
$t=0$.  Finally, of those eigenvalues which vanish to order
$n-1$, the difference between the number of eigenvalues
whose $n^{\hbox{\eightrm{th}}}$ derivative is positive and
the number whose $n^{\hbox{\eightrm{th}}}$ derivative is
negative is given by the signature of the form
$Q_n:\calg^{\ell-1}_n\times\calg^{\ell-1}_n\ra\bbbC$ which
is independent of the metric on $M$.

In this section we will show that for a path
$d_t:\Omega^{*}(M;E_{\bbbC})\ra\Omega^{*}(M;E_{\bbbC
})$ of flat connections, the functions $\delta_n$ defined in
section 2 can be identified with certain higher Massey
products which are well-known to be invariants of the
homotopy type of $M$.

The deRham complex $\{\Omega^{*}(M;E_{\bbbC}),d_0\}$ is
a differential graded Lie algebra and it is in this context
that we will now work.

\definition{}{A {\bfit graded Lie algebra} (\gla) is a
sequence $\{L^n\}^{\infty}_{n=0}$ of complex vector spaces
together with a bilinear operation $[\;\;,\;\;]:L^m\times
L^n\ra L^{m+n}$ satisfying:

\item{(1)}$[x,y]+(-1)^{mn}[y,x]=0$

\item{(2)}$[x,[y,z]]=[[x,y],z]+(-1)^{mn}[y,[x,z]]$

\noindent where we are assuming $x\in L^m$ and $y\in L^n$.
(2) is just the graded version of the Jacobi identity.  We
write $L=\opluslim\limits^{\infty}_{n=0}L^n$, and extend
$[\;\;,\;\;]$ linearly over $L$.}

A {\bfit differential graded Lie algebra} (\dgla) is a \gla\
equipped with a differential
$$d:L^m\ra L^{m+1}$$
for all $m$ which satisfies $d[x,y]=[dx,y]+(-1)^{|x|}[x,dy]$.
 ($|x|$ is
defined by $x\in L^{|x|}$.)

Assume $(L,d)$ is a \dgla.  For each $a\in L^1$, we define
$d_a:L^m\ra L^{m+1}$ by $d_ax=dx+[a,x]$.  The Jacobi
identity implies that $(L,d_a)$ is also a \dgla.  We call
the \dgla\ $(L,d)$ {\bfit flat} if $d\circ d=0$.  If $(L,d)$
is flat then $(L,d_a)$ is flat if and only if $da+{1\over
2}[a,a]=0$.

Note that if $A$ is a flat connection, then the deRham
complex $(\Omega^{*}(M;E_{\bbbC}),d_A)$ introduced in
section 5 is a flat \dgla.  It is a well-known fact
that the exterior derivatives $d_{A^{\prime}}$, coming from
other connections $A\rq$ on $P$, are precisely those
differentials of the form $d_{A^{\prime}}x=d_Ax+[a,x]$
where $a\in\Omega^{1}(M;E_{\bbbC})$.  We say the form
$a$ is {\bfit flat} with respect to a fixed flat connection
$A$ if $d_Aa+{1\over 2}[a,a]=0$.

Let $(L,d)$ be a flat \dgla\ throughout the rest of the
section.  Let
$L[[t]]\equiv\opluslim\limits^{\infty}_{n=0}L^n[[t]]$ denote
the space of formal power series with coefficients in $L$.
By defining $[xt^n,yt^m]=[x,y]t^{n+m}$ and $d(xt^n)=d(x)t^n$
we give $(L[[t]],d)$ the structure of a flat \dgla.  If
$a(t)=\sum\limits^{\infty}_{i=1}a_it^i\in L^1[[t]]$ then we
may define $d_{a(t)}:L[[t]]\ra L[[t]]$ as before by
$d_{a(t)}x(t)=dx(t)+[a(t),x(t)]$.  Note that $a(t)$ is flat
if and only if the following sequence of equations holds:
$$\openup3\jot\eqalign{
da_1&=0\cr
da_2&=-{1\over 2}[a_1,a_1]\cr
da_3&=-{1\over 2}([a_1,a_2]+[a_2,a_1])\cr
&\vdots\cr
da_n&=-{1\over 2}\sum\limits^{n-1}_{i=1}[a_i,a_{n-i}]\cr}$$

Given a flat $a(t)\in L^1[[t]]$, with vanishing constant
term we now define a sequence of relations $T_n$ on $L$.
These relations will coincide, in the case where $a(t)$ is
convergent for $t\in(-\ve,\ve)$, with the relations defined
in section 2 for the operator $d_tx=dx+[a(t),x]$.

Define $uT_nv$ if and only if there exist $u_1,\ldots,u_{n-
1}\in L$ such that
$$d_{a(t)}(u+u_1t+\ldots+u_{n-1}t^{n-1})=vt^n+w(t)t^{n+1}$$
for some $w(t)\in L[[t]]$.

If we define $\calz_n$, $\calb_n$,
$\calg_n=\calz_n/\calb_n$, and $\delta_n:\calg_n\ra\calg_n$
just as in section 2, (with $\delta_n$ induced by $n!T_n$)
then the proof of Theorem 2.1 implies that
$(\calg_n,\delta_n)$ is a sequence of well-defined chain
complexes with the chain groups of each complex equal to the
homology of the one before, i.e., $\calg_{n+1}=\ker
\delta_n/\Im\delta_n$.   Next, we will show that in the
context of a flat \dgla\ the operator $\delta_n$ can be
expressed in terms of Massey products as defined by Retakh in
[Re].

We begin by reviewing Retakh's definition of higher Massey
products in a flat \dgla.  We define a {\bfit multi-index}
$\title I$ with $|I|=n\geq 1$ to be a set of integers
$I=\{i_1,\ldots,i_n\}$ such that $1\leq i_1<i_2<\ldots<i_n$.
A {\bfit submulti-index of} $\title I$ is a non-empty subset
of $I$.  A {\bfit proper pair} of submulti-indices of $I$ is
an ordered pair $(J,K)$ of submulti-indices of $I$, where we
write
$$J=\{j_1<\ldots<j_r\}\qquad{\hbox{\rm{and}}}\qquad
K=\{k_1<\ldots<k_s\}$$
which satisfy

\item{(1)}$J\cap K=\emptyset$

\item{(2)}$J\cup K=I$

\item{(3)}$J\not=\emptyset\not=K$

\item{(4)}$j_1<k_1$ (i.e., $j_1=i_1$)

\noindent The set of all proper pairs of submulti-indices of
$I$ is denoted $PP(I)$.

For each homogeneous element $x\in L$, let $\barx=(-
1)^{|x|+1}x$.  Let $x_1,\ldots,x_n$ be an ordered $n$-
tuple of homogeneous cocycles in $L$.  A {\bfit Massey
system} $\calm$ for $(x_1,\ldots,x_n)$ consists of a
homogeneous element $u_I\in L$ for each multi-index
$I\subset\{1,2,\ldots,n\}$ satisfying the following two
properties:

\item{(1)}$u_{\{i\}}=x_i$ for each $i\in\{1,\ldots,n\}$

\item{(2)}For $|I|>1$,
$$du_I=\sum\limits_{(J,K)\in PP(I)}(-
1)^{\ve(J,K)}[\baru_J,u_K]$$
where $\ve(J,K)=\sum\limits_{(j,k)\in J\times
K:k<j}(|x_j|+1)(|x_k|+1)$.

Given a Massey system $\calm=\{u_I\}$ for
$\{x_1,\ldots,x_n\}$, the corresponding Massey product is
defined by
$$\tilQ_{\calm}(x_1,\ldots,x_n)=\sum\limits_{(J,K)\in
PP\{1,\ldots,n\}}(-1)^{\ve(J,K)}[\baru_J,u_K]$$

In practice, to construct a Massey system one constructs the
elements $u_I$ inductively according to $|I|$.  To begin, if
$1\leq i_1<i_2\leq n$, one selects a homogeneous element
$u_{i_1,i_2}$ such that
$$du_{i_1,i_2}=[\barx_{i_1},x_{i_2}].$$
To accomplish this step, it must be true that
$[x_{i_1},x_{i_2}]=0$ in $H^*(L)$.  Of course, this
condition only determines $u_{i_1,i_2}$ up to addition of a
$1$-cocycle.  Having defined all the $u_{i_1,i_2}$s, one
next must choose $u_{i_1,i_2,i_3}$ (for $1\leq
i_1<i_2<i_3\leq n$) satisfying
$$du_{i_1,i_2,i_3}=[\baru_1,u_{23}]+[\baru_{12},u_3]+(-
1)^{(|x_2|+1)(|x_3|+1)}[\baru_{13},u_2],$$
et cetera.

It is a routine exercise (using the Jacobi identity) to show
that once one has defined the elements $u_J$ for all $|J|<m$
then, if $|I|=m$, $\sum\limits_{(J,K)\in PP(I)}(-
1)^{\ve(J,K)}[\baru_J,u_K]$ is a cocycle of degree
$(\sum\limits_{i\in I}|x_i|)-|I|+1$.  Hence in order to
define $u_I$, this cocycle must be a coboundary.

As a special case of the above ``routine exercise'' we have
the following proposition.

\proposition{6.1}{If $\calm$ is a Massey system for
$\{x_1,\ldots,x_n\}$, then $\tilQ_{\calm}(x_1,\ldots,x_n)$
is a cocycle of degree $(\sum\limits_{1\leq i\leq n}|x_i|)-
n+1$.}

\definition{}{Let $Q_{\calm}(x_1,\ldots,x_n)$ denote the
cohomology class in $H^*(L)$ represented by
$\tilQ_{\calm}(x_1,\ldots,x_n)$.}

In general $Q_{\calm}(x_1,\ldots,x_n)$ depends on which
Massey system $\calm$ was chosen.
Note that a Massey system for $(x_1,\ldots,x_n)$ contains Massey
systems for all proper subsets of $(x_1,\ldots,x_n)$ for which
the corresponding Massey products vanish.  Thus the Massey
product cannot be defined for an arbitrary $n$-tuple of
cocycles.  The second order Massey product is, however,
well-defined for every pair $(x_1,x_2)$ of cocycles; it is
simply the ``cup product'' $[\barx_1,x_2]$.

\proposition{6.2}{Let the power series
$a(t)=\sum\limits^{\infty}_{i=1}a_it^i$ in $L^1[[t]]$ be
flat in the sense that $d_{a(t)}\circ d_{a(t)}=0$, as
discussed earlier.  For all multi-indices $I$, define
$u_I=(-1)^{n+1}n!a_n$, where $n=|I|$.  Then

\item{(1)}$u_{\{i\}}=a_1$ for all $i\in\bbbZ_+$

\noindent and

\item{(2)}$du_I=\sum\limits_{(J,K)\in PP(I)}(-
1)^{\ve(J,K)}[\baru_J,u_K]$ for all $i\subseteq\bbbZ_+$.

}

In words, this proposition states that the elements $u_I$,
as defined in the proposition, contain a Massey system for
the ordered $n$-tuple $x_1=x_2=\ldots=x_n=a_1$ for every
$n$, and all of the corresponding Massey products vanish.

\proof Set $c_n=(-1)^{n+1}n!a_n$ for each $n\geq 1$.  Assume
$|I|=n$.  Compute $du_I=dc_n=(-1)^{n+1}n!da_n=(-1)^{n+1}n!(-
{1\over 2})\sum\limits^{n-1}_{i=1}[a_i,a_{n-i}]={1\over
2}\sum\limits^{n-1}_{i=1}{n\choose i}[c_i,c_{n-
i}]=\sum\limits_{(J,K)\in PP(I)}[u_J,u_K]$, where the last
step follows by an easy combinatorial argument.  Because
$u_I\in L^1$ for all $I$, it follows that $\baru_I=u_I$ and
$(-1)^{\ve(J,K)}=1$, so the proposition follows.\proofend

Continue to assume that
$a(t)=\sum\limits^{\infty}_{i=1}a_it^i$ is a flat series in
$L^1[[t]]$.  Let $u$ be a homogeneous cocycle (of any
degree) and note that $a_1$ is also a cocycle since $a(t)$
is flat.  We will now show that $uT_nv$ if and only if
$\tilQ_{\calm}(u,{\underbrace{a_1,\ldots,a_1}_n})=(-1)^{n+1}n!v$ for
some Massey system $\calm$ which is $a(t)$-compatible in the
following sense:

\definition{}{Suppose $n\geq 1$, $a(t)$ is a flat series in
$L^1[[t]]$, $u$ is a homogeneous cocycle, $\calm=\{u_I\}$
is a Massey system for
$(u,{\underbrace{a_1,\ldots,a_1}_n})$.  Then we say $\calm$
is $a(t)$-compatible if and only if

\item{(1)}if $1\not\in I$, then $u_I=(-1)^{n+1}n!a_n$ where
$n=|I|$

\noindent and

\item{(2)}if $1\in I_1\cap I_2$ and $|I_1|=|I_2|$, then
$u_{I_1}=u_{I_2}$.

}

More informally, this definition says that for $\calm$ to be
$a(t)$-compatible:  (1) $\calm$ must make use of the Massey
systems contained in $a(t)$ to the largest extent possible;
and (2) $\calm$ must have as much symmetry as possible.

\proposition{6.3}{Let $d_t=d_{a(t)}:L\ra L$, where $a(t)\in
L^1[[t]]$ is flat.  Then $uT_nv$ if and only if
$\tilQ_{\calm}(u,{\underbrace{a_1,\ldots,a_1}_n})=(-
1)^{n+1}n!v$ for some $a(t)$-compatible Massey system
$\calm$.}

\proof We will prove the ($\Rightarrow$) direction.  The
proof of ($\Leftarrow$) is just a reversal of the argument
and will be omitted.  Recall that $uT_nv$ means there exists
$u_1,\ldots,u_{n-1}\in L$ such that
$$d_t(u+u_1t+\ldots+u_{n-1}t^{n-1})=vt^n+w(t)t^{n+1}.$$
This implies that
$$\openup3\jot\eqalign{
&du=0\cr
&du_1+[a_1,u]=0\cr
&\vdots\cr
&du_{n-1}+[a_1,u_{n-2}]+\ldots+[a_{n-1},u]=0.\cr}$$
It follows that if we set $u_I=(-1)^{n+1}n!a_n$ (where
$n=|I|$) if $1\not\in I$ and $u_I=(-1)^nn!u_n$ (where
$n=|I|-1$) if $1\in I$ then $\calm=\{u_I\}_I$ s an
$a(t)$-compatible Massey system and
$\tilQ(u,{\underbrace{a_1,\ldots,a_1}_n})=(-1)^{n+1}n!v$.
To verify this last fact, note that
$$\tilQ_{\calm}(u,a_1,\ldots,a_1)=\sum\limits_{(J,K)\in
PP\{1,\ldots,n+1\}}(-1)^{\ve(J,K)}[\baru_J,u_K].$$
$\ve(J,K)=1$ because $|a_1|=1$, and
$$\baru_J=(-
1)^{|u_J|+1}u_J=(-1)^{|u_J|+1}u_J=(-1)^{|u|+1}u_J,$$ so
$$\eqalign{\tilQ_{\calm}(u,a_1,\ldots,a_1)&=(-
1)^{|u|+1}\sum\limits_{(J,K)\in
PP\{1,\ldots,n+1\}}[u_J,u_K]\cr
&=(-1)^{|u|+1}\sum\limits^{n-
1}_{i=0}{n\choose i}[(-1)^i(i!)u_i,(-1)^{n-i+1}(n-i)!a_{n-
1}\cr
&=n!(-1)^{n+1}\sum\limits^{n-1}_{i=0}[a_{n-i},u_i]\cr
&=(-1)^{n+1}n!v.\cr}$$
  We omit the routine verification that $\calm$
is an $a(t)$-compatible Massey system.  This completes the
proof of the ($\Rightarrow$) direction.\proofend

Recall from section 2 that $\tau_n:\calg_n\ra\calg_n$ is
just the well-defined map induced by the relation $T_n$ on
$V$, and $\delta_n=n!\tau_n$.  Hence for $[u]\in\calg^i_n$,
we may write
$$\delta_n[u]=(-
1)^{n+1}[\tilQ_{\calm}(u,{\underbrace{a_1,\ldots,a_1}_n})]
\in\calg^{i+1}_n$$
Note the following consequence of what we have proved:  the
normally ill-defined Massey product
$$u\mapsto Q_{\calm}(u,{\underbrace{a_1,\ldots,a_1}_n})$$
becomes well-defined when we insist that $\calm$ be
$a(t)$-compatible, and when we think of the domain and range of the
function as being the subquotients $\calg^i_n$ of
$H^i(M;E_{\bbbC})$ and $\calg^{i+1}_n$ of
$H^{i+1}(M;E_{\bbbC})$.

\titlea{7}{General bundles}

In the previous two sections we chose  the bundle $E$ to be the
complexification of the adjoint bundle associated to some principal bundle.

This was done  for simplicity and we indicate now what modifications need
to be made to handle the case when $E$ is a general  vector bundle with a
path of flat connections.

Let $P\to M$ be a principal $G$ bundle for some compact lie group $G$ with
Lie
algebra $\calg$, and let
$\rho:G\to U(W)$ be a unitary representation of $G$ onto a
hermitian vector
space $W$, defining a vector bundle $E=P\times_\rho W$.  Let $adP$ denote
the
adjoint bundle of lie algebras, $adP=P\times_{ad} \calg$.
 An analytic path
of
flat connections on $P$  defines an analytic path of operators
$d_t^E:\Omega^*(M;E)\to \Omega^*(M;E)$, as well as an analytic path of
operators
$d_t^{ad}:\Omega^*(M;adP)\to \Omega^*(M;adP)$.  Then
$a(t)=d_t^{ad}-d_0^{ad}$
is an analytic path in $\Omega^1(M;adP)$ and $d_t^E=d_0^E+ \rho_*(a(t))$
where
$\rho_*:\calg\to gl(W)$ is the differential of $\rho$.
For simplicity of
notation we assume that $\calg$ acts on $W$ on the {\it right}.

The differential of $\rho$ endows $\Omega^*(M;E)$ with the
 structure  of a
module over the DGLA $\Omega^*(M;adP)$.
Thus to extend the definitions of
Massey products we need a ``hybrid'' generalization
 which takes this action
into account.  This is done as follows.

If $x_1\in \Omega^*(M;E)$  is a $d_0^E$-cocycle and
$x_i\in\Omega^*(M;adP), \ i=2,\cdots n$ are $d_0^{ad}$-cocycles, we define
Massey systems as  before, with the caveats:
\item{1.}   If $1\in I$, then $u_I\in \Omega^*(M;E)$, if $1\not\in I$,
then
$u_I\in \Omega^*(M;adP)$.
\item{2.} If $1\not\in J\cup K$, then $[\bar{u}_J,u_K]$ has the same
meaning as
before.  If $1\in J$, then $[\bar{u}_J,u_K]$ means $\bar{u}_J \rho_*(u_K)$,
i.e. the action of $u_K$ on $\bar{u}_J$.
\item{3.} In the equation
$$du_I=\sum\limits_{(J,K)\in PP(I)}(-
1)^{\ve(J,K)}[\baru_J,u_K],$$
one should interpret the operator $d$ as $d_0^E$ or $d_0^{ad}$
according to
whether or not $1\in I$.

With this definition, one can define Massey products
$${\calq}_{\calm}(e,y_1,\cdots,y_{n})\in H^*(M;E)_{A(0)}$$
of cocycles $e\in \Omega^*(M;E)$ and $y_i\in \Omega^*(M;adP)$.

 The definition
of $a(t)$-compatible Massey systems for $\calq(u,a_1,\cdots, a_1)$  is
similar: one assumes that $u_I=(-1)^{n+1}n!a_n$ if $1\not\in I$, and  if
$1\in I_1\cap I_2$ and $|I_1|=|I_2|$, then $u_{I_1}=u_{I_2}$.  Then the
theorem
of the previous section continues to hold with these definitions.

Finally, to define the bilinear forms $Q_n$, one need only use the ``dot''
product coming from a $G$-invariant fiber metric on the bundle $E$. With
these
definitions, all the results of Section 5 and 6 continue to hold for $E$ a
general vector bundle.

\vfill\eject
\centerline{\bf REFERENCES}
\vskip5ex

\item{[Fa]}M. Farber, {\sl Singularities of the analytic torsion}, 1994,
preprint.

\item{[FL]}M. Farber and J. Levine,  {\sl Deformations of the
Atiyah-Patodi-Singer eta-invariant}, 1993, preprint.

\item{[Ka]}T. Kato, {\sl Perturbation theory for linear
operators}, Grundlehren der Mathematischen Wissenschaft 132,
Springer Verlag, Berlin, 1966.

\item{[KK1]}P. Kirk and E. Klassen,  {\sl Computing spectral flow via cup
products}, 1993, to appear in J. Diff. Geom.

\item{[KK2]}P. Kirk and E. Klassen,  {\sl The spectral flow of the
odd signature operator on manifolds with boundary}, 1993, preprint.

\item{[Mc]}McCleary, {\sl A user's guide to spectral sequences}, Publish
or Perish, Wilmington, 1985.

\item{[Re]}Retakh, {\sl Massey operations in Lie superalgebras and
deformations of complexly analytic algebras}, Functional Analysis
and its Applications, Apr. 1979, p. 319--322. (Russian Original: Vol. 12,
No.
4, Oct.- Dec. 1978.)

\baselineskip=12pt
\vskip\baselineskip
\noindent Paul A. Kirk

\noindent Department of Mathematics

\noindent Indiana University

\noindent Bloomington, IN   47405

\vskip\baselineskip

\noindent Eric P. Klassen

\noindent Department of Mathematics

\noindent Florida State University

\noindent Tallahassee, FL  32306-3027

\vfill\end